\magnification=1200
\input amstex
\input amssym.def
\documentstyle{amsppt}
\parindent=0pt
\font\fabs=cmr8 scaled\magstep0

\define\la{\langle}
\define\ra{\rangle}
\define\bs{\backslash}
\define\Z{{\Cal Z}_n}
\redefine\ze{\zeta}
\define\z+{{\Bbb Z}_+}
\define\Hn{{\Cal H}_n}
\define\HTn{{\widetilde{\Cal H}}_n}
\define\HTnn{{\widetilde{\Cal H}}_{n-1}}
\define\Hnn{{\Cal H}_{n-1}}
\define\A_n{{\Cal A}_q(n)}
\define\U_n{{\Cal U}_q(n)}
\define\UU{{\Cal U}_q(n-1)}
\define\ep{\varepsilon}
\define\zw{z^{\lambda}w^{\mu}}
\redefine\H{Hopf $\ast$-algebra\ }
\define\wi{\widetilde}
\define\ov{\overline}
\define\hb{\hfill\break}
\define\ho{\natural}
\redefine\a{\alpha}
\redefine\C{{\Bbb C}}
\redefine\b{\beta}
\redefine\g{\gamma}
\redefine\d{\delta}
\define\pr{{^{\prime}}}
\define\ot{\otimes}
\define\sq{$\hfill \square$}

\headline={\ifodd\pageno\rightheadline\else\leftheadline\fi}
\def\rightheadline{\tenrm\hfil Addition Formula for $q$-Disk Polynomials 
\hfil\folio}
\def\leftheadline{\tenrm\folio\hfil Paul G.A. Floris \hfil}
\voffset=2\baselineskip

\topmatter
\title Addition Formula for $q$-Disk Polynomials\endtitle
\author Paul G.A. Floris\endauthor
\affil Leiden University \\ Department of Mathematics and Computer Science \\
P.O. Box 9512 \\ 2300 RA Leiden\\ The Netherlands\\ e-mail: 
floris\@wi.leidenuniv.nl\\ fax: +31 - 71 - 27 69 85 \endaffil
\endtopmatter

{\fabs {{\sl Keywords and phrases} : CQG algebras, quantum unitary group, 
quantized universal enveloping algebra, spherical elements, 
$q$-disk polynomials, addition formula.}}
\vskip0.5cm
{\fabs {{\sl 1991 Mathematics Subject Classification} : 16W30, 17B37, 33D45,
33D80, 81R50.}}
\bigskip

{\bf 0. Introduction.}
\medskip
Over the past few years quantum groups have shown to be powerful tools in
the study of $q$-hypergeometric functions.
They enabled proofs of identities which would have been hard to guess
without quantum group theoretic motivation.
As an example one can mention the papers [Ko2], [Koo5], [Koo6]
and [NM]. See also [Ko3], [Koo3] and [N1] for surveys on the connection between
quantum groups and basic hypergeometric functions and references therein.
\smallskip
The purpose of this paper is to present an addition formula
for so-called $q$-disk polynomials, using some quantum group theory.
This result is a $q$-analogue of a result which was proved around
1970 by ${\breve{\text S}}$apiro [S] and Koornwinder [Koo1,2] independently.
They considered the homogeneous space $U(n)/U(n-1)$, were $U(n)$
denotes the group of unitary transformations of the vector space $\C^n$, and
identified the corresponding zonal spherical functions as disk polynomials,
in the terminology of [Koo2].
These are orthogonal polynomials in two variables whose 
orthogonality measure is supported by the closed unit disk in the complex
plane, and which can be expressed in terms of Jacobi polynomials $P_n^{(\a,\b)}
(x)$
for certain integer values of the parameters $\a$ and $\b$. The associated spherical
funtions were shown to be also expressible in terms of disk polynomials, and
from this an addition formula was proved for (positive) integer values 
of $\a$ and $\b$. By an easy argument this identity was then extended 
to all complex values of $\a$ (and from this Koornwinder could even derive 
an addition formula for general Jacobi polynomials). And in fact the line of 
arguing and the results in this paper will be very similar to (part of) the
ones in [Koo2].
\smallskip
$q$-Disk polynomials (quantum disk polynomials in the terminology of [Koo4])
are polynomials in two non-commuting variables which are expressed by
means of little $q$-Jacobi polynomials, and that can be understood as a 
$q$-analogue in non-commuting variables of disk polynomials.
They appeared in [NYM] where Noumi, Yamada and Mimachi studied
a quantum analogue of $U(n)/U(n-1)$, or, rather, the coordinate ring of
this quantum homogeneous space. They investigated its 
${\Cal U}_q(\frak g
\frak l(n))$-module structure, where ${\Cal U}_q(\frak g\frak l(n))$ denotes
the quantized universal enveloping algebra corresponding to the Lie algebra
$\frak g\frak l(n)$, using the somewhat more abstract theory of highest weight
representations, and ended with identifying the zonal spherical functions
as $q$-disk polynomials.
\smallskip
This paper is organized as follows. In chapter one we recall some of the
facts on compact quantum groups and quantum homogeneous spaces.
In the second chapter we introduce a $q$-deformation $\Z$ of the algebra of
polynomials on ${\C }^n$ and study its structure as a ${\Cal U}_q(\frak 
g\frak l(n))$-module. Chapter three then deals with $\wi{\Z}$, the $q$-deformed 
algebra of polynomials on the sphere in ${\C }^n$. This algebra is the same
as the algebra $A(K\backslash G)$ of [NYM, section 4.1]. In the first two
sections we introduce invariant integration on $\wi{\Z}$ and describe its
irreducible decomposition as a ${\Cal U}_q(\frak g\frak l(n))$-module.
In section three we recover the zonal spherical elements as $q$-disk 
polynomials, using the invariant integral. Section four deals with the
irreducible decomposition as a ${\Cal U}_q(\frak g\frak l(n-1))$-module
and the associated spherical elements. These turn out to be expressable
through $q$-disk polynomials as well. Finally, in the fifth section we
prove the addition formula for the $q$-disk polynomials.\hb
In view of the previous paragraph, part of this paper (up till section 3.3) 
does not contain any
really new results, but merely presents an alternative and somewhat more
explicit approach to the problems addressed in [NYM].
The new results are contained in sections 3.4 and 3.5, where we recover
the
associated spherical functions and 
establish the addition formula. It should be noted that the
definition of $q$-disk polynomials as polynomials in two non-commuting 
variables accounts for the fact that the addition formula we end up with 
is an identity in several non-commuting variables as well.
\medskip
We end by introducing some notation which will be frequently used.\hb
Let $0<q<1$. Write $\z+$ for the set of non-negative integers: 
$$
\z+ = \{0,1,2,\hdots\}.
$$
For $l,m\in \z+$, $0<q<1$, $\b\in\C$, $a$ an indeterminate and $\lambda = 
(\lambda_1, \hdots, \lambda_n)\in \z+^n$ we define the following symbols:
$$
\eqalign{
l\wedge m &= \text{min} (l,m) = {1\over 2}(l+m - \vert l-m\vert)\cr
[m]_q &= {{1 - q^m}\over {1-q}} \cr
(a;q)_{\b} &= {{(a;q)_{\infty}}\over {(aq^{\b};q)_{\infty}}}\cr
(a;q)_{\infty} &=\prod_{k=0}^{\infty} (1-q^ka)\cr
\vert\lambda\vert &= \lambda_1 + \hdots +\lambda_n.\cr}
$$
Furthermore we recall Jackson's $q$-integral for a function $f$
$$
\int_0^c f(x)d_q x = c(1-q) \sum_{k=0}^{\infty} f(cq^k) q^k,
$$
and the little $q$-Jacobi polynomials
$$
p_m(x;a,b;q) = {_2}\varphi_1 \bigl[ {{q^{-m}\ ,\ abq^{m+1}}\atop {aq}}; 
q,qx\bigr] = \sum_{k=0}^m {{(q^{-m};q)_k (abq^{m+1};q)_k}\over {(aq;q)_k 
(q;q)_k}} (qx)^k.
$$
If $0<aq<1$ and $bq<1$ they satisfy the orthogonality
$$
\sum_{k=0}^\infty {{(bq;q)_k}\over{(q;q)_k}}(aq)^k \big( p_lp_m\bigr)
(q^k;a,b;q) =\d_{l,m} {{(q,bq;q)_l (aq)^l (1-abq) (abq^2;q)_\infty}\over
{(aq,abq;q)_l (1-abq^{2l+1}) (aq;q)_\infty}}.
$$
In particular, if we let
$$
P_m ^{(\a,\b)} (x;q) = p_m(x;q^{\a},q^{\b};q)\qquad\quad (\a,\b>-1)
$$
then the orthogonality reads
$$
\eqalign{
\int_0^1 P_l^{(\a,\b)}(x;q) &P_m^{(\a,\b)}(x;q) x^{\a} 
{{(q x;q)_{\infty}}\over {(q^{\b+1}x;q)_{\infty}}} d_{q} x =\cr
&\d_{lm}{{(1-q)q^{m(\a+1)}}\over {1-q^{\a+\b+2m+1}}} {{(q;q)_m (q;q)_{\b+m}}\over 
{(q^{\a+1};q)_m (q^{\a+1};q)_{\b+m}}}.\cr}
$$
\medskip
Through the little $q$-Jacobi polynomials one can define the $q$-disk
polynomials as follows (cf. [Koo4]):\hb
Suppose we are given a complex unital $\ast$-algebra $\Cal Z$ generated by 
the elements $z$ and $z^{\ast}$, subject to the relation
$$
z^{\ast} z= q^2 zz^{\ast} + 1-q^2
$$
and with $\ast$-structure $(z)^{\ast}=z^{\ast}$. Then the $q$-disk polynomials
$R_{l,m}^{(\a)} (z,z^{\ast};q)$, with $\a>-1$ and $l,m\in \z+$, are defined 
to be:
$$
R_{l,m}^{(\a)}(z,z^{\ast};q) =\left\{\aligned z^{l-m} P_m^{(\a,l-m)}
(1-zz^{\ast};q) \qquad (l\geq m)\\
P_l^{(\a,m-l)}(1-zz^{\ast};q) (z^{\ast})^{m-l}\qquad (l\leq m).
\endaligned\right.
$$
\smallskip
Throughout this paper we will keep $0<q<1$ fixed.
\medskip
{\bf Acknowledgements }: The author is very grateful to Prof. Tom H.
Koornwinder for his valuable comments and suggestions.

\vfill\eject

{\bf 1. CQG algebras}.
\medskip
This first chapter will be devoted to establishing some notations and to 
recall 
briefly some facts about compact quantum groups. Our language will be that of 
[DK]. See also [Koo7] for a more extensive treatment of this topic.\hb
The ground field for all vector spaces under consideration will be the 
field $\C$ of complex numbers.
\medskip
{\bf 1.1. Definitions and generalities}.
\medskip
Let $A$ be a Hopf algebra with comultiplication $\Delta: A \to A\ot A$, 
counit $\ep : A\to \C$ and antipode $S : A\to A$ (for the theory of Hopf 
algebras we refer the reader to [Sw]).\hb
A {\it right corepresentation} of $A$ is a pair $(V,\pi)$ of a complex
 vector space $V$ and a linear map $\pi : V\to V\ot A$, satisfying 
$$
(\pi\ot id)\pi = (id\ot\Delta)\pi,\ \ (id\ot \ep)\pi = id.\leqno(1.1.1)
$$
We also say that $V$ is a right comodule for $A$. If $V$ is finite dimensional 
with basis $\{ e_i\}_{i=1}^N$ and if we write $\pi(e_i) =\sum_{k=1}^N e_k\ot 
\pi_{ki}$, then this is equivalent to saying that the $\pi_{ij}$ satisfy 
$\Delta(\pi_{ij}) =\sum_{k=1}^N \pi_{ik}\ot \pi_{kj}$ and $\ep(\pi_{ij})= 
\delta_{ij}$. The elements $\pi_{ij}$ of $A$ are called the matrix 
coefficients of 
this corepresentation. Furthermore, an element $v$ of $V$ is said to be 
{\it (right) $A$-invariant} if $\pi(v) = v\ot 1$.\hb
A particular example of a right corepresentation is given when $V$ is a 
subspace of $A$ such that $\Delta(V)\subset V\ot A$ and $\pi=\Delta$. In 
this case we will also call $V$ a (right) coideal of $A$. Similar 
definitions can be given for left corepresentations, but we won't be needing 
them here.\hb
A Hopf algebra $A$ is called a {\it Hopf $\ast$-algebra} if there exists 
an anti-linear involution $\ast : A\to A$ which turns $A$ into a 
$\ast$-algebra and such that $\Delta$ and $\ep$ are $\ast$-homomorphisms. 
It is not difficult to show that in this case the antipode $S$ is invertible 
and satisfies $S\circ\ast\circ S\circ\ast = id$. \hb
Now suppose $A$ is a \H. A right corepresentation $\pi : V\to V\ot A$ is 
called {\it unitarizable} if there exists a hermitean inner product 
$\la .,.\ra$ on $V$ such that
$$
\la \pi(v), \pi(w)\ra = \la v,w\ra 1_A
$$
for all $v,w\in V$. Here we extended the inner product $\la .,.\ra$ to a map 
from $V\ot A$ to $A$ by letting
$$
\la v\ot a,w\ot b\ra = \la v,w\ra b^{\ast}a.
$$
If $V$ is endowed with this inner product, then the corepresentation is 
called {\it unitary} and the inner product is called {\it $A$-invariant}. 
If in this setting $V$ is finite dimensional with orthonormal basis 
$\{e_i\}$ with respect to a given inner product, then the following are 
equivalent:
\parindent=8pt
\roster
\item "(i)" $(V,\pi)$ is a unitary corepresentation
\item "(ii)" $\sum_{k=1}^N \pi_{ki}^{\ast}\pi_{kj}^{} = \delta_{ij}1_A$
\item "(iii)" $S(\pi_{ij}) = \pi_{ji}^{\ast}$
\item "(iv)" $\sum_{k=1}^N \pi_{ik}^{}\pi_{jk}^{\ast} =\delta_{ij}1_A$.
\endroster
\parindent=0pt
\hb
Let us write $\Sigma=\Sigma(A)$ for the set of equivalence classes of finite 
dimensional irreducible unitary corepresentations of the \H $A$. Furthermore, 
for $\pi =(\pi_{ij})_{i,j=1}^N \in\Sigma$ let us put $A_{\pi}= \text{span} 
\{\pi_{ij}\}_{i,j=1}^N$ and $A_{\pi(r)} = \text{span} \{\pi_{rj}\}_{j=1}^N$ 
$(1\leq j\leq N)$. Then one can prove that the $\pi_{ij}$ $(\pi\in\Sigma)$ 
are linearly independent and $\sum_{\pi\in\Sigma} A_{\pi}$ is a direct sum 
(see e.g. [Koo7]; actually this is even true in the more general situation of 
coalgebras).
\smallskip
{\bf Definition} : A \H $A$ is called a {\it CQG algebra} if $A$ is spanned 
by the matrix coefficients of all its finite dimensional (irreducible) 
unitary corepresentations, i.e. if $A=\sum_{\pi\in\Sigma} A_{\pi}$.
\smallskip
The direct sum decomposition $A=\sum_{\pi\in\Sigma} A_{\pi}$ is often referred 
to as the Peter-Weyl decomposition. We should also note that with respect 
to $\Delta$ we have the following irreducible decomposition of $A$ as a 
right comodule
$$
A= \bigoplus_{\pi\in \Sigma} \bigoplus_{r=1}^{d_{\pi}} A_{\pi(r)}.
$$
A particular class of CQG algebras is given by the finitely generated ones, 
called CMQG algebras. A CQG algebra $A$ is a CMQG algebra if and only if there 
exists a single finite dimensional unitary corepresentation $t =(t_{ij})$ 
of $A$ such that $A$, as an algebra, is generated by the matrix elements 
$t_{ij}$. This corepresentation is usually referred to as the {\it fundamental 
(or natural) corepresentation}.\hb
As an example of a CMQG algebra we will meet the algebra ${\Cal A}_q(n) = 
{\Cal A}_q(U(n))$ of regular functions on the quantum unitary group.
\medskip
\proclaim\nofrills{\bf Theorem 1.1.1}: Let $A$ be a CQG algebra. Then there 
exists a unique linear functional $h:A\to\C$, called normalized Haar 
functional, that satisfies
\parindent=8pt
\roster
\item "(i)" $h(1_A) =1$
\item "(ii)" $(h\ot id)\Delta (a) = h(a)1_A = (id\ot h)\Delta (a)$
\item "(iii)" $h(a^{\ast}a)>0$ if $a\neq 0$.
\endroster
\parindent=0pt
\endproclaim
\medskip
\proclaim\nofrills{\bf Theorem 1.1.2} : If $A$ is a CQG algebra, then every 
finite dimensional 
corepresentation is unitarizable. Moreover, for all $\pi\in\Sigma$ there 
exists a positive definite 
matrix $G^{\pi}$ such that if $\pi,\rho\in\Sigma$ then
$$
\eqalign{
h(\pi_{ij}^{}\rho_{kl}^{\ast}) &= \delta_{\pi\rho}\delta_{ik} {{G^{\pi}_{jl}}
\over 
{\text{tr }G^{\pi}}}\cr
h(\pi_{ij}^{\ast}\rho_{kl}^{}) &= \delta_{\pi\rho}\delta_{jl} 
{{(G^{\pi})^{-1}_{ki}}\over 
{\text{tr } (G^{\pi})^{-1}}}.\cr}
$$
\endproclaim
\smallskip
Proofs of these results can be found in [Koo7]. The last result is referred 
to as Schur orthogonality.
\medskip

{\bf 1.2 Transitive coactions}.
\medskip
In this section we will adopt the terminology of [D] section 4.1.\hb
Let a CQG algebra $A$ and a $\ast$-algebra $Z$ be given. Furthermore let us 
assume that there exists a 
{\it right $\ast$-coaction} $\delta : Z\to Z\ot A$ of $A$ on $Z$. This means 
that $\delta$ is a 
corepresentation in the sense of section 1.1 and moreover a homomorphism of 
$\ast$-algebras. Suppose that 
this coaction is {\it transitive}, i.e. that there exists an injective 
$\ast$-algebra homomorphism 
$\psi : Z\to A$ which intertwines the coactions $\delta$ on $Z$ and $\Delta$ 
on $A$. Then we know 
from [D], Theorem 4.1.5, that $Z$ possesses a normalized, positive definite 
{\it $A$-invariant functional} 
$h_Z : Z\to \C$, meaning 
that $h_Z$ satisfies $(h_Z\ot id)\delta (z) = h_Z(z) 1_A$ for all $z\in Z$, 
that $h_Z(1_Z)=1$ and
$h_Z(z^{\ast}z)>0$ if $z\neq 0$. Moreover, the coaction $\d$ 
is unitary with respect to the inner product on $Z$ given by
$$
\la z,w\ra = h(w^{\ast}z).\tag 1.2.1
$$
Upon identifying $Z$ with $\psi (Z)$, we may assume that $Z$ is a 
$\ast$-subalgebra and right coideal 
of $A$, and that $\delta = \Delta$ and $h_Z = h$.
\medskip
Assume there exists a \H $C$ with unital Hopf $\ast$-algebra epimorphism 
$\theta : A\to C$ and such 
that for all $z\in Z$ there holds
$$
\theta (z) = \ep (z) 1_C.\tag 1.2.2
$$
Using this mapping $\pi$, we can define a coaction $\delta_C$ from $C$ on $Z$ by putting $\delta_C  = 
(id\ot \theta)\Delta$.
\medskip
\proclaim\nofrills {\bf Proposition 1.2.1}: Suppose $V$ is a subcomodule of 
$Z$ of finite dimension 
$N\geq 1$, so $\Delta : V\to V\ot A$. Then $V$ contains a non-zero $C$-invariant vector. That is, there 
exists an element $\ze_0\neq 0$ in $V$ with the property $\delta_C (\ze_0) = 
\ze_0\otimes 1_C$.
\endproclaim
\smallskip
{\sl Proof}: Let $\{ e_i\}_{i=1}^N$ be an orthonormal basis for $V$ with 
respect to the inner product 
$(1.2.1)$. As before we write $\Delta(e_i) =\sum_{k=1}^N e_k\ot \pi_{ki}$. 
Unitarity of $\Delta$ will 
imply that $\Delta\boxtimes\Delta := (id\ot id\ot m_A) (id\ot \sigma_{23}\ot 
id)(\Delta\ot \Delta) : V\ot 
V\to V\ot V\ot A$ is 
also a unitary coaction with respect to the inner product $\la v_1\ot w_1, 
v_2\ot w_2\ra = \la v_1, 
v_2\ra \la w_1, w_2\ra$. Here $m_A : A\ot A\to A$ denotes multiplication in
$A$.
Furthermore, we have
$$
\Delta\circ m_A = (m_A\ot id)\circ\Delta\boxtimes\Delta.\tag 1.2.3
$$
Now consider the element $\zeta := \sum_{i=1}^N e_i\ot e_i^{\ast}$ in $Z\ot Z$. 
It satisfies
$$
\eqalign{
\Delta\boxtimes\Delta (\zeta) &= \sum_{i=1}^N \Delta\boxtimes\Delta (e_i\ot 
e_i^{\ast})\cr
&= \sum_{i,j,k=1}^N e_j\ot e_k^{\ast}\ot \pi_{ji}^{}\pi_{ki}^{\ast}\cr
&= \sum_{j,k=1}^N  \delta_{jk} e_j\ot e_k^{\ast}\ot 1_A\cr
&= \zeta\ot 1_A \cr}\tag 1.2.4
$$
i.e. $\zeta$ is $A$-invariant. If we put ${\tilde {\zeta}} := m (\zeta) 
=\sum_{i=1}^N e_ie_i^{\ast}\in Z$, 
then it follows from $(1.2.3)$ and $(1.2.4)$ that
$$
\Delta ({\tilde {\zeta}}) = \Delta\circ m_A(\zeta) = 
(m_A\ot id)\Delta\boxtimes\Delta(\zeta) = {\tilde 
{\zeta}}\ot 1_A,
$$
and hence ${\tilde {\zeta}} = (id \ot h)\Delta({\tilde {\zeta}})$. By the 
invariance of the Haar 
functional we get ${\tilde {\zeta}} = h({\tilde {\zeta}}) 1_A$. From this 
we obtain
$$
\eqalign{
\ep ({\tilde {\zeta}})
&= h({\tilde {\zeta}})\cr
&= \sum_{i=1}^N h(e_ie_i^{\ast})\cr
&=  \sum_{i=1}^N \la e_i^{\ast},e_i^{\ast}\ra\cr
&>0. \cr}\tag 1.2.5  
$$
Finally, put $\ze_0:= (id\ot\ep)(\zeta) = \sum _{i=1}^N \ep (e_i^{\ast}) e_i$. 
Then obviously $\ze_0\in 
V$. Also $\ze_0\neq 0$, since $\ep (\ze_0) = \ep ({\tilde {\zeta}}) >0$.\hb
The last thing to show is that $\ze_0$ is $C$-invariant. Writing $m_C$ for 
the multiplication in $C$, 
one has
$$
\eqalign{
\d_C (\ze_0) &= (id\ot\theta)\Delta(\ze_0)\cr
&= (id\ot m_C)\bigl( (id\ot\theta)\Delta(\ze_0)\ot 1_C\bigr)\cr
&= (id\ot m_C)\bigl( (id\ot\theta)(\Delta\ot\ep)(\ze)\ot 1_C\bigr)\cr
&= (id\ot m_C)\bigl( (id\ot\theta\ot id)(\Delta\ot (\ep\ot\theta)\Delta)(\ze)\cr
&= (id\ot m_C)\bigl( (id\ot\theta\ot\ep\ot\theta)(\Delta\ot\Delta)(\ze)\cr
&= (id\ot m_C)\bigl( (id\ot\ep\ot\theta\ot\theta) (id\ot\sigma_{23}\ot id)(\Delta\ot\Delta)(\ze)\cr
&= (id\ot\ep\ot\theta) (id\ot id\ot m_A)(id\ot\sigma_{23}\ot id)(\Delta\ot\Delta)(\ze)\cr
&= (id\ot\ep\ot\theta) (\Delta\boxtimes\Delta)(\ze)\cr
&= (id\ot\ep\ot\theta)(\ze\ot 1_A)\cr
&= (id\ot\ep)(\ze)\ot 1_C\cr
&= \ze_0\ot 1_C.\cr}
$$
This proves the proposition.\sq 
\medskip
{\bf Remark}: Note that the element $\ze$ does not depend on the choice of 
the orthonormal basis of $V$. 
Suppose the pair $(A,C)$ forms a so-called {\it quantum Gel'fand pair}, so 
that 
in each $A_{\pi}$ there 
is an at most one-dimensional subspace of $C$-invariant elements.
If the comodule structure on $V$ is unitary and irreducible, and we consider 
$V$ as a subspace of $A$ as 
before, then $V$ can be realized as some $A_{\pi(r)}$, say $A_{\pi(1)}$, for 
certain $\pi\in\Sigma$ (cf. 
section 1.1). The invariant subspace  will then be spanned by $\pi_{11}$. 
From Theorem 1.1.2 we know 
that the elements $e_i = \bigl( {{\text{tr}(G^{\pi})^{-1}}\over 
{((G^{\pi})^{-1})_{11}}}\bigl)^{1\over 2} 
\pi_{1i}$ form an orthonormal basis of $V$. Thus, again by Theorem 1.1.2
$$
\ep(\ze_0) = \sum_{i=1}^N h(e_ie_i^{\ast}) = {{\text{tr}(G^{\pi})^{-1}}\over 
{((G^{\pi})^{-1})_{11}}} = 
{1\over {h(\pi_{11}^{\ast}\pi_{11})}}.
$$
This remark will be of importance in section 3.5.
\medskip
{\bf Remark}: The element $\ze$ plays a role similar to the kernel function 
on a compact homogeneous 
space: let 
$T_{\ze} :Z\to Z$ be the mapping 
$T_{\ze} (z) = (id\ot h_Z)\bigl( \ze .(1\ot z)\bigr)$. Then it is easily seen 
that $T_{\ze} (v) =v$ 
for all $v\in V$ and that $T_{\ze}(w)=0$ if $w$ is orthogonal to $V$. Indeed, 
if $v=\sum_{j=1}^N \a_j 
e_j$, then $T_{\ze}(v) = \sum_{i,j=1}^N \a_j e_i 
h_Z(e_i^{\ast}e_j) = \sum_{i,j=1}^N \a_j e_i \la e_i,e_j\ra = 
\sum_{j=1}^N \a_j e_j =v$. In the same way
it follows that $T_{\ze}(w)=0$ if $\la w ,V\ra=0$.
\bigskip
We will end this section with a small lemma, needed further on.
\smallskip
\proclaim\nofrills {\bf Lemma 1.2.2}: Suppose $A,C$ are two CQG-algebras with 
normalized Haar functionals 
$h_A, h_C$ respectively, and such that there exists a surjective unital Hopf 
$\ast$-algebra homomorphism 
$\pi:A\to C$ ($C$ is a 'subgroup' of $A$). Suppose furthermore that $Z$ is a 
right comodule for $A$
under $\delta: Z\to Z\ot A$, and assume that $\d$ is unitary with respect to 
a given 
inner product $\la . , .\ra$.\hb
Let $V\subset Z$ be a subcomodule for the corepresentation $\d_C = (id\ot\pi)\d$, and let a 
$C$-invariant element $\phi$ in $Z$ be given.\hb
Now if $\phi$ is orthogonal with respect to $\la . , . \ra$ to all $C$-invariant elements in $V$, then 
$\phi$ is orthogonal to the whole of $V$.
\endproclaim
\medskip
{\sl Proof }: We know that $\delta_C(\phi) = \phi\ot 1_C$. Consider the linear 
map $T_C : Z\to Z$ 
defined by $T_C = (id\ot h_C)\delta_C$. Since $\phi$ is $C$-invariant, one has 
$T_C(\phi) = \phi$. 
Furthermore, if $\psi\in V$ is arbitary, then $T_C(\psi)$ is $C$-invariant 
and in $V$. Now note 
that the $A$-invariant inner product $\la .,.\ra$ is also $C$-invariant: 
$\la \delta_C(\phi), 
\delta_C (\psi)\ra = \la \phi,\psi\ra 1_C$ for all $\phi, \psi \in Z$.
Consequently, if $\phi\in Z$ is $C$-invariant and orthogonal to all 
$C$-invariant elements in $V$, 
and if $\psi\in V$ is arbitrary then
$$
\eqalign{
\la \phi, \psi\ra &= h_C(\la \phi,\psi\ra 1_C)\cr
&= h_C(\la \delta_C(\phi), \delta_C(\psi)\ra)\cr
&= h_C( \la \phi\ot 1_C, \delta_C (\psi)\ra )\cr
&= \la \phi, T_C(\psi)\ra\cr
&=0\cr}
$$
since $T_C(\psi)$ is $C$-invariant. This proves the lemma.\sq
\bigskip

{\bf 2. The quantized algebra of polynomials on $\C^n$}.
\medskip
In this chapter we will define the algebra $\Z$, which is a $q$-deformation 
of the involutive algebra 
of polynomials on $\C^n$. On this algebra we will define a $\ast$-action of 
$\U_n$, the quantized 
universal enveloping algebra of the unitary group $U(n)$, as the 'differential' 
of a certain 
$\ast$-coaction of $\A_n$, the quantized algebra of regular functions on $U(n)$. 

\vfill\eject

{\bf 2.1. Definition and structure of $\Z$}.
\smallskip
Define the $\ast$-algebra $\Z$ as $\Z = \Bbb C\la z_1,\hdots, z_n, 
w_1,\hdots,w_n\ra / I$, where $I$ 
is the two-sided ideal generated by all elements
$$
\eqalign{
z_iz_j - qz_jz_i\qquad (1\leq i<j\leq n)\cr
w_jw_i - qw_iw_j\qquad (1\leq i<j\leq n)\cr
w_iz_j - qz_jw_i\qquad (1\leq i,j\leq n,i\neq j)\cr
w_iz_i - z_iw_i - (1-q^2)\sum_{k<i} z_kw_k\quad (1\leq i\leq n).\cr}
$$
The $\ast$-structure on $\Z$ is given by $z_i^{\ast} = 
w_i \ (1\leq 1\leq n).$\hb
In other words, $\Z$ is the complex involutive algebra generated by the 
elements $z_i,w_i\ 
(1\leq i\leq n)$, subject to the relations
$$
\leqalignno{
z_iz_j &= qz_jz_i\qquad (1\leq i<j\leq n)&(2.1.1)\cr
w_jw_i &= qw_iw_j\qquad (1\leq i<j\leq n)&(2.1.2)\cr
w_iz_j &= qz_jw_i\qquad  (1\leq i,j\leq n,i\neq j)&(2.1.3)\cr
w_iz_i &= z_iw_i + (1-q^2)\sum_{k<i} z_kw_k\quad (1\leq i\leq n)&(2.1.4)\cr}
$$
and with involution $\ast:\Z\to\Z,\ z_i^{\ast}=w_i$. For $q=1$ this algebra can 
be interpreted as the 
commutative involutive algebra of polynomials in the $n$ coordinates 
$z_1,\hdots, z_n$ on $\C^n$ and 
their complex conjugates.
\medskip
\proclaim\nofrills {\bf Proposition 2.1.1}: $\Z$ has as a $\Bbb C$-linear basis 
the collection
$$
\{ z^{\lambda}w^{\mu} := z_1^{\lambda_1}\hdots z_n^{\lambda_n}w_n^{\mu_n}\hdots 
w_1^{\mu_1} | 
\lambda,\mu\in \z+^n \}.
$$
\endproclaim
{\sl Proof}: This is proved with the aid of Bergman's Diamond Lemma (see [B]).
\sq
\medskip
{\bf Remark}: Similarly one can prove that $\Z$ has a $\Bbb C$-basis 
given by $\{ w^{\mu}z^{\lambda} 
:= w_1^{\mu_1}\hdots w_n^{\mu_n}z_n^{\lambda_n}\hdots z_1^{\lambda_1} | 
\lambda,\mu\in 
\z+^n \}.$ This also follows from Proposition 2.1.1 by applying the 
algebra isomorphism 
which interchanges $z_i$ and $w_i$ $(1\leq i\leq n)$ and sends $q$ to $q^{-1}$.
\medskip
For $i=1,\hdots,n$ define the following elements in $\Z$:
$$
Q_i = \sum_{k=1}^i z_kw_k.
$$
Then we have the easily verified identities:
\smallskip
\proclaim\nofrills {\bf Lemma 2.1.2}: For $1\leq i,k\leq n$ and $m\in \z+$ 
there holds
$$
\alignat 2
Q_i^{\ast} &= Q_i & \qquad Q_iQ_k &= Q_k Q_i\\
z_k w_k &= Q_k - Q_{k-1}  & w_kz_k &= Q_k - q^2 Q_{k-1}\\
z_k Q_i &= q^{-2} Q_i z_k \qquad\  \text{and} & w_k Q_i &= q^2 Q_i w_k 
\qquad {\text if}\ \  k>i\\
z_k Q_i &= Q_i z_k\qquad\qquad \text{and} & w_k Q_i &= Q_i w_k \qquad\quad 
{\text if}\ \ k\leq i\\ 
z_k^mw_k^m &= Q_k^m ({{Q_{k-1}}\over {Q_k}} ; q^{-2})_m &\quad w_k^m z_k^m 
&= Q_k^m (q^2 {{Q_{k-1}}\over 
{Q_k}} ; q^2)_m.
\endalignat
$$
\endproclaim
As to the center of $\Z$ one has
\proclaim\nofrills {\bf Proposition 2.1.3}: The element $Q_n$ is central 
in $\Z$, and for the center Cent$(\Z)$ of $\Z$ we have
$$
\text{Cent} (\Z) = \Bbb C [Q_n].\tag 2.1.5
$$
\endproclaim
{\sl Proof}: It is easy to check that $Q_n$ commutes with all the $z_i$ and 
$w_i$, and therefore 
that $\C [Q_n]\subset \text{Cent}(\Z)$. To prove the second part of the 
proposition we use an 
induction argument. Introduce a total ordering on the basis of $\Z$ as follows: 
to the basis element 
$\zw$ associate the sequence $\{\lambda,\mu\} = (|\lambda |+|\mu |,
 \lambda_n,\hdots,\lambda_1,
\mu_1,\hdots,\mu_n)$. Now we declare $\zw \succeq z^{\rho}w^{\sigma}$ if 
and only if $\{\lambda,\mu\} 
\geq \{\rho,\sigma\}$ with respect to the lexicographic ordering of 
elements in 
$\z+^{2n+1}$.\hb
Let $\phi\in\text{Cent}(\Z)$ be non-zero. Write
$$
\phi=c\zw + \sum_i c_i z^{\lambda^{(i)}}w^{\mu^{(i)}}
$$
where $c\neq 0$ and $\zw \succ z^{\lambda^{(i)}} w^{\mu^{(i)}}$ for all $i$. 
Denote by $\ho (\phi)$ 
the highest order part of $\phi$ with respect to the ordering as above, so 
$\ho (\phi) =c \zw$. Being central, $\phi$ must satisfy 
$\phi z_i = z_i\phi$ for all $1\leq i\leq n$. So in particular $\ho (\phi z_i) 
=\ho (z_i \phi)$. But
$$
\eqalign{
\ho (\phi z_i) &= cq^{\mu_1+\hdots +\mu_{i-1}+\mu_{i+1}+\hdots +\mu_n -
(\lambda_{i+1}+\hdots +\lambda_n)} 
z^{\lambda +\ep_i}w^{\mu}\cr 
\ho (z_i \phi) &= cq^{-(\lambda_1 +\hdots +\lambda_{i-1})} z^{\lambda +
\ep_i}w^{\mu}.\cr}
$$
where $\lambda +\ep_i = (\lambda_1,\hdots,\lambda_i +1,\hdots, \lambda_n)$.
Since $q$ is not a root of unity this yields
$$
\lambda_1 +\hdots +\lambda_{i-1} =-(\mu_1+\hdots +\mu_{i-1}+\mu_{i+1}+\hdots 
+\mu_n) +\lambda_{i+1}+
\hdots +\lambda_n\tag 2.1.6
$$
for all $1\leq i\leq n$.
Taking $i=n$ we get $0=\lambda_1 +\hdots +\lambda_{n-1} + \mu_1 +\hdots 
+\mu_{n-1}$, hence
$$
\lambda_j =\mu_j =0 \qquad (1\leq j\leq n-1).\tag 2.1.7
$$
But combining the case $i=n-1$ of $(2.1.6)$ with $(2.1.7)$ gives 
$\lambda_n=\mu_n$. Thus
$$
\phi=cz_n^{\lambda_n}w_n^{\lambda_n} + \sum_i c_i z^{\lambda^{(i)}}w^{\mu^{(i)}}.
$$
On the other hand, $\ho (cQ_n^{\lambda_n}) = cz_n^{\lambda_n}w_n^{\lambda_n}$. 
So letting $\psi = 
\phi -cQ_n^{\lambda_n}$ we see that $\ho (\psi) \prec \ho (\phi)$. 
Applying the 
induction hypothesis, 
$\psi\in \C [Q_n]$ and hence $\phi = \psi + cQ_n^{\lambda_n}\in\C [Q_n]$. 
This completes the proof.\sq
\medskip
{\bf Remark }: What we actually constantly use is the following fact: 
if $\phi\in\Z$ is a monomial which
contains $z_i$ and $w_i$ exactly $\lambda_i$ respectively $\mu_i$ times 
($1\leq i\leq n)$, then $\phi$ can
be written such that $\ho (\phi) =c\zw$ with $c\neq 0$ and $\lambda = (\lambda_1,
\hdots, \lambda_n),\ \mu=(\mu_1,\hdots,\mu_n)$.

\vfill\eject

{\bf 2.2. A ${\Cal U}_q(\frak g\frak l (n))$-module structure on $\Z$.}
\medskip
Define the \H ${\Cal A}_q(n) = {\Cal A}_q (U(n))$ as the free unital complex algebra on the 
generators $t_{ij} (1\leq i,j \leq n), det_q^{-1}$, subject to the relations
$$
R^+T_2T_1 = T_1T_2R^+\ \ ,\ \ det_q^{-1}det_q = 1 = det_q det_q^{-1}.
$$
Here $R^+$ is the R-matrix for root system $A_{n-1}$ given as
$$
R^+ = \sum_{i,j=1}^n q^{\delta_{ij}} E_{ii}\ot E_{jj} + (q-q^{-1}) \sum_{1\leq i<j\leq n} 
E_{ij}\ot E_{ji}
$$
and $det_q = det_q(T)$ is the so-called quantum determinant which generates the center of 
${\Cal A}_q(n)$. It is given by
$$
det_q = \sum_{\sigma\in S_n} (-q)^{l(\sigma)} t_{\sigma(1) 1}\hdots t_{\sigma(n) n}
$$
in which we let $S_n$ denote the group of permutations of the set $\{ 1,\hdots, n\}$ and 
$l(\sigma)$ the smallest number of transpositions $(i\ i+1)$ in which $\sigma\in S_n$ can be decomposed. 
For the Hopf $\ast$-structure on ${\Cal A}_q(n)$ and further details we refer the reader to [N2]
or [NYM]. 
$\A_n$ is an example of a CMQG-algebra (cf. section 1.1).\hb
To this \H corresponds a dual \H, the quantized universal enveloping algebra $\U_n = {\Cal U}_q 
(\frak g\frak l (n))$ with generators $q^h\ (h\in P^{\ast}),\ e_k,\ f_k\ (k=1,\hdots,n-1)$. 
Here $P= \sum_{i=1}^n {\Bbb Z}\ep_i$ is the weight lattice for $\frak g\frak l(n)$ and 
$P^{\ast}= \text{Hom}(P,\Bbb Z)$ is identified with $P$ via the non-degenerate bilinear 
form $\la \ep_i,\ep_j\ra =\delta_{ij}$. Again we refer the reader to [N2] for the exact 
definition and structure maps on this \H. Note that in particular we define the  comultiplication 
$\Delta$ on the generators $q^h,\ e_k$ and $f_k$ as
$$
\eqalign{
\Delta (q^h) &= q^h\ot q^h\cr
\Delta (e_k) &= q^{\ep_k -\ep_{k+1}}\otimes e_k + e_k \otimes 1\cr
\Delta (f_k) &= 1\otimes f_k + f_k\otimes q^{-({\ep_k -\ep_{k+1}})}.\cr}\tag 2.2.1
$$
The \H duality between ${\Cal A}_q(n)$ and $\U_n$ is defined on the generators by
$$
\eqalign{
\la q^h, t_{ij}\ra &= \delta_{ij} q^{\la h,\ep_i\ra}\cr
\la e_k, t_{ij}\ra &= \delta_{ki}\delta_{k+1,j}\cr
\la f_k, t_{ij}\ra &= \delta_{k+1,i}\delta_{kj}.\cr}\tag 2.2.2
$$
Now we have a right $\ast$-coaction of ${\Cal A}_q(n)$ on the generators of $\Z$ as follows:
$$
\eqalign{
\delta : \Z &\to \Z\otimes {{\Cal A}_q(n)}\cr
z_i &\mapsto \sum_{k=1}^n z_k \otimes t_{ki}\cr
w_i &\mapsto \sum_{k=1}^n w_k \otimes t_{ki}^{\ast}\cr}\tag 2.2.3
$$
linearly extended in both factors. 
\proclaim\nofrills {\bf Lemma 2.2.1}: The map $\delta$ defined as in $(2.2.3)$ extends to a 
$\ast$-algebra homomorphism on $\Z$ and satisfies (1.1.1).
\endproclaim
\smallskip
{\sl Proof} : Extending $\delta$ as a $\ast$-algebra homomorphism, we only need to verify 
that it respects the relations (2.1.1-4). One can readily check that this is true if and only 
if we have the following identities in ${\Cal A}_q(n)$:
$$
\alignat 2
t_{ki}t_{kj} &= qt_{kj}t_{ki} &\qquad &(i<j)\\
t_{ki}t_{lj}-t_{lj}t_{ki} &= qt_{kj}t_{li} - q^{-1}t_{li}t_{kj}&\qquad 
&(i<j,\ k<l)\\
t_{ki}^{\ast}t_{lj}^{} &= t_{lj}^{}t_{ki}^{\ast}&\qquad &(i\neq j,\ k\neq l)\\
t_{ki}^{\ast}t_{kj}^{} &= qt_{kj}^{}t_{ki}^{\ast} -(1-q^2)\sum_{l>k} 
t_{li}^{\ast}t_{lj}^{}&\qquad 
&(i\neq j)\\
qt_{ki}^{\ast}t_{li}^{} &= t_{li}^{}t_{ki}^{\ast} +(1-q^2)\sum_{j<i} 
t_{lj}^{}t_{kj}^{\ast}&\qquad 
&(i\neq j)\\
t_{ki}^{\ast}t_{ki}^{} + (1-q^2)\sum_{l>k}t_{li}^{\ast}t_{li}^{} &= t_{ki}^{}t_{ki}^{\ast} + 
(1-q^2)\sum_{j<i} t_{kj}^{}t_{kj}^{\ast}.&&
\endalignat
$$
And this is indeed the case (see for example [Ko1], section 2). It is now 
immediate that $\d$ satisfies $(1.1.1)$.\sq
\smallskip
"Differentiating" this right coaction gives a left $\ast$-action of the quantized universal 
enveloping algebra $\U_n$:
$$
X.\phi = (id\otimes X)\delta (\phi)\qquad (X\in \U_n, \phi\in \Z).\tag 2.2.4
$$
Here the dot . denotes the action, and the $X$ in the right-hand side should be understood as a 
linear functional on ${\Cal A}_q(n)$ induced by the Hopf $\ast$-algebra pairing $\la .,. \ra: 
\U_n\times {\Cal A}_q(n)\to \Bbb C$ as given by $(2.2.2)$. 
\smallskip
\proclaim\nofrills {\bf Lemma 2.2.2}:  $(2.2.4)$ defines an algebra action of $\U_n$ on $\Z$, 
and satisfies $X.(\phi\psi) = \sum _{(X)} \bigl(X_{(1)}.\phi\bigr) \bigl(X_{(2)}.\psi\bigr)$. 
Moreover $X.\phi^{\ast} = (S(X)^{\ast}.\phi)^{\ast}$ for all $X\in \U_n$ and all $\phi,\psi\in 
\Z$. Here we use the symbolic notation $\Delta (X) = \sum_{(X)} X_{(1)}\otimes X_{(2)}$.
\endproclaim
\smallskip
Note that if $\phi\in\Z$ is ${\Cal A}_q(n)$-invariant, this will imply by definition that 
$X.\phi =\ep(X)\phi$ for all $X\in \U_n$. In this case we will say that $\phi$ is $\U_n$-invariant.
\smallskip
\proclaim\nofrills {\bf Lemma 2.2.3} : The element $Q_n$ is ${\Cal A}_q(n)$-invariant and hence 
also $\U_n$-invariant.
\endproclaim
{\sl Proof} : This follows immediately from $(2.2.3)$, the fact that $\delta$ is an algebra 
homomorphism and the relation $\sum_{k=1}^n t_{ik}^{}t_{jk}^{\ast}=\d_{ij}1$
(see [Ko1, (2.12)]).\sq
\medskip
The corresponding $\ast$-action of $\U_n$ is given in the following proposition:
\smallskip
\proclaim\nofrills {\bf Proposition 2.2.4}: For all $h\in P^{\ast},\ 1\leq k\leq n-1$ and all 
$\lambda,\mu\in \z+^n$ there holds
$$
\eqalign{
q^h. z^{\lambda}w^{\mu} &= q^{\la h,\lambda -\mu\ra} z^{\lambda}w^{\mu}\cr
f_k. z^{\lambda}w^{\mu} &= -q^{\mu_k +1} [\mu_{k+1}]_{q^{-2}} z^{\lambda}w^{\mu -\ep_{k+1} 
+\ep_k} + q^{\lambda_{k+1} +\mu_k -\mu_{k+1}} [\lambda_k]_{q^{-2}} z^{\lambda -\ep_k +
\ep_{k+1}}w^{\mu}\cr
e_k.z^{\lambda}w^{\mu} &= -q^{-1}q^{\mu_{k+1} +\lambda_k -\lambda_{k+1}} [\mu_k]_{q^{-2}} 
z^{\lambda}w^{\mu +\ep_{k+1}-\ep_k} + q^{\lambda_k} [\lambda_{k+1}]_{q^{-2}} z^{\lambda +\ep_k 
-\ep_{k+1}}w^{\mu}.\cr}
$$
Here $\lambda +\ep_i = (\lambda_1,\hdots ,\lambda_i +1,\hdots, \lambda_n)$.
\endproclaim
\smallskip
{\sl Proof}: First prove, using $(2.2.1-4)$ and Lemma 2.2.2, that for all $h\in P^{\ast},\ 
1\leq k\leq n-1$, $1\leq i\leq n$ and $m\in \z+$ we have
$$
\eqalign{
q^h.z_i^m &= q^{m\la h,\ep_i\ra}z_i^m\cr
q^h.w_i^m &= q^{-m\la h,\ep_i\ra}w_i^m\cr
f_k.z_i^m &= \delta_{ik} [m]_{q^{-2}} z_i^{m-1}z_{i+1}\cr
f_k.w_i^m &= -q \delta_{i,k+1} [m]_{q^{-2}} w_i^{m-1}w_{i-1}\cr
e_k.z_i^m &= \delta_{i,k+1} [m]_{q^{-2}} z_{i-1}z_i^{m-1}\cr
e_k.w_i^m &= -q^{-1}\delta_{ik} [m]_{q^{-2}} w_{i+1}w_i^{m-1}\cr}
$$
by using induction with respect to $m$. From this one gets that
$$
\eqalign{
q^h.z^{\lambda} &= q^{\la h,\lambda\ra} z^{\lambda}\cr
q^h.w^{\mu} &= q^{-\la h,\mu\ra} w^{\mu}\cr
f_k.z^{\lambda} &= q^{\lambda_{k+1}} [\lambda_k]_{q^{-2}} z^{\lambda -\ep_k +\ep_{k+1}}\cr
f_k.w^{\mu} &= -q^{\mu_k +1} [\mu_{k+1}]_{q^{-2}} w^{\mu -\ep_{k+1} + \ep_k}\cr
e_k.z^{\lambda} &= q^{\lambda_k} [\lambda_{k+1}]_{q^{-2}} z^{\lambda +\ep_k -\ep_{k+1}}\cr
e_k.w^{\mu} &= -q^{-1}q^{\mu_{k+1}} [\mu_k]_{q^{-2}} w^{\mu +\ep_{k+1} - \ep_k}.\cr}
$$
Combination of this with $(2.2.1)$ will conclude the proof.\sq
\smallskip
{\bf Remark}: One can rewrite this last result a little when using that for all $m\in \z+$ one 
has the identity $[m]_{q^{-2}} = q^{-2(m-1)} [m]_{q^2}$.
It is actually this left-action of $\U_n$ that we will consider rather than the $\A_n$-coaction, 
since it is the more easy to handle.
\medskip
We have the following decomposition of $\Z$:
$$
\Z = \bigoplus_{l,m\in \z+} \Z (l,m)
$$
where $\Z (l,m)$ is the subspace of $\Z$ spanned by all elements which are homogeneous of degree 
$l$ in the $z_k$ and homogeneous of degree $m$ in the $w_k$ . Note that it makes sense to speak 
of homogeneous elements, since the relations $(2.1.1-4)$ are homogeneous.
Also note that $0$ is homogeneous of any degree $(l,m)$. It follows from Proposition 2.1.1 that $\Z 
(l,m)$ has a linear basis consisting of all elements $\zw$ with the property that $|\lambda |=l, 
|\mu |=m$. Hence the $\Z (l,m)\ (l,m\in \z+)$ are finite dimensional. It is easy to check that 
they are subcomodules of $\Z$ under $\delta$ (and hence submodules for the action of $\U_n$). 
Note also that
$$
Q_n \Z (l-1,m-1)\subset \Z (l,m)
$$
and
$$
\Z (l,m)^{\ast} = \Z (m,l).
$$
\medskip
\proclaim\nofrills {\bf Proposition 2.2.5}: Suppose $\phi\in \Z (l,m)$ is $\U_n$-invariant. Then:
$$
\eqalign{
&(i)\ \text{if } l\neq m,\ \text{then } \phi=0\cr
&(ii)\ \text{if } l = m,\ \text{then } \phi=c Q_n^l\qquad (c\in {\Bbb C}).\cr}
$$
\endproclaim
{\sl Proof}: We already know from Lemma 2.2.3 that $Q_n$ is $\U_n$-invariant. Suppose 
$\phi\in \Z (l,m)$ is $\U_n$-invariant, and write $\phi=c\zw + \sum_i c_i z^{\lambda^{(i)}}w^{\mu^{(i)}}$ 
where $c\neq 0$ and for all $i$ there holds $\zw \succ z^{\lambda^{(i)}}w^{\mu^{(i)}}$ with respect 
to the total ordering on momomials in $\Z$ (cf. Proposition 2.1.3). As before 
let $\ho (\phi)$ be the 
highest order part of $\phi$ (so $\ho (\phi)= c\zw$). \hb
Since $\ep (q^h)=1$ for all $h\in P^{\ast}$, we must have $\ho (q^{\ep_i}.\phi) 
= \ho (\phi)$. But, 
by Proposition 2.2.4, we see that $\ho (q^{\ep_i}.\phi) = cq^{\lambda_i - \mu_i}\zw$. Hence this 
implies that $\lambda_i=\mu_i$ for all $1\leq i\leq n$, and therefore $l= |\lambda | =|\mu | =m$ 
if $\phi\neq 0$. This proves part $(i)$.\hb
Next we prove $(ii)$ using induction with respect to $\{\lambda,\mu\}$. Suppose $l=m$ and suppose 
that $\phi = c\zw + \psi$ is $\U_n$-invariant, but 
not equal to a constant times $Q_n^l$. Here $\psi$ is a sum of terms which are smaller that $\zw$ 
with respect to the total ordering on monomials. As in $(i)$ we conclude that $\lambda =\mu$, 
so $\phi= cz^{\lambda}w^{\lambda} + \psi$. We may assume $cz^{\lambda}w^{\lambda}\neq z_n^lw_n^l$ 
(if not, then consider $\phi - cQ_n^l$, which is $\U_n$-invariant). This means there exists a $k_0$, 
with $1\leq k_0\leq n-1$, such that $\lambda_{k_0}\neq 0$. But since $\ep (e_k) =0$ $(1\leq 
k\leq n-1)$, we must have that $\ho (e_k.\phi) =0$ for $k=1,\hdots,n-1$. Hence, by Proposition 2.2.4,
we get
$$
-cq^{-1}q^{{\lambda}_k}[\lambda_k]_{q^{-2}} z^{\lambda}w^{\lambda + \ep_{k+1} -\ep_k} = 0\qquad 
(k=1,\hdots,n-1)
$$
i.e.
$$
\lambda_k = 0 \qquad (k=1,\hdots,n-1).
$$
So in particular 
$\lambda_{k_0} =0$, yielding a contradiction. This proves part $(ii)$, and thus the proposition.\sq
\medskip
\proclaim\nofrills {\bf Lemma 2.2.6}: For all $1\leq i\leq n$ and each $m\in \z+$ we have
$$
w_i^mz_i = q^{2m} z_iw_i^m + (1-q^{2m}) w_i^{m-1} Q_i.
$$
\endproclaim
\smallskip
{\sl Proof}: Note that $w_iz_i = q^2 z_iw_i + (1-q^2) Q_i$. Now proceed by induction with 
respect to $m$.\sq
\medskip
\proclaim\nofrills {\bf Corollary 2.2.7}: For all $1\leq i\leq n$ and each $m\in\z+$ one has
$$
(z_iw_i)^m = \sum_{k=0}^m c_k z_i^kw_i^k Q_i^{m-k}
$$
for certain coefficients $c_k\in\C [q]$.
\endproclaim
\medskip
\proclaim\nofrills {\bf Corollary 2.2.8}: For all $1\leq i\leq n$ and all $m,p\in\z+$ one has
$$
w_i^mz_i^p = \sum_{k=0}^{m\wedge p} d_k z_i^{p-k}w_i^{m-k}Q_i^k
$$
for certain $d_k\in\C [q]$.
\endproclaim
\medskip
{\bf Remark}: Actually one has that $c_k, d_k\in {\Bbb C} [q^2]$ for all $k$.
\medskip
Suppose we are given the two algebras $\Z$ and ${\Cal Z}_s$ with generators $z_i,\ w_i\ 
(1\leq i\leq n)$ and $z_i\pr,\ w_i\pr\ (1\leq i\leq s)$ respectively and suppose $s<n$. 
Then we have the following canonical embedding
$$
\iota^{(s,n)} : {\Cal Z}_s\to \Z
$$
which sends the generators $z_i\pr, w_i\pr$ of ${\Cal Z}_s$ to the first $s$ pairs of generators 
$z_i,\ w_i\ (1\leq i\leq s)$ of $\Z$, and the restriction map
$$
\rho^{(n,s)} : \Z\to {\Cal Z}_s
$$
which puts $z_i$ and $w_i$ equal to zero for $i=1,\hdots, n-s$ and maps $z_i,\ w_i$
to $z_{i-n+s}\pr,\ w_{i-n+s}\pr$ respectively for $i=n-s+1,\hdots, n$. 
Both maps are $\ast$-algebra homomorphisms. So in particular we can view 
${\Cal Z}_{n-1} (l,m)$ as lying in $\Z (l,m)$, using $\iota^{(n-1,n)}$.\hb
Furthermore, observe that for $1\leq p\leq n-1$ we have a natural embedding ${\Cal U}_q(n-p)
\hookrightarrow \U_n$ by identifying ${\Cal U}_q(n-p)$ with the subalgebra of $\U_n$ generated 
by the elements $q^{\ep_i}\ (1\leq i\leq n-p), e_k, f_k\ (1\leq k\leq n-p-1)$. In this way it 
is possible to speak of ${\Cal U}_q(n-p)$-invariant elements in $\Z$.
\medskip 
\proclaim\nofrills {\bf Proposition 2.2.9}: Suppose $\phi\in \Z (l,m)$ is ${\Cal U}_q(n-1)$-invariant. 
Then $\phi$ is of the form
$$
\phi = \sum_{j=0}^{l\wedge m} c_j z_n^{l-j}w_n^{m-j} Q_n^j.\tag 2.2.5
$$
Conversely, any $\phi$ of this form is $\UU$-invariant. Hence the dimension of $\UU$-invariant 
elements in  $\Z (l,m)$ equals $l\wedge m + 1$.
\endproclaim
{\sl Proof}: Using Proposition 2.1.1 and the commutation relations for the $z_i, w_i$, we see 
that $\phi$ can be uniquely written as
$$
\phi = \sum_{i=0}^l\sum_{j=0}^m z_n^{l-i}w_n^{m-j}p_{ij}(z_1,\hdots,z_{n-1},w_1,\hdots,w_{n-1})
$$
for certain $p_{ij}\in {\Cal Z}_{n-1} (i,j)$. Then the action of ${\Cal U}_q(n-1)$ is on the 
elements $p_{ij}$ (cf. Proposition 2.2.4). By Proposition 2.2.5 we get that $p_{ij} =\delta_{ij} 
d_jQ_{n-1}^j$, and thus $\phi = \sum_{j=0}^{l\wedge m} d_j z_n^{l-j}w_n^{m-j} Q_{n-1}^j$, which 
already yields the stated dimension. Since $Q_n$ is central, we can write
$Q_{n-1}^j = (Q_n - z_nw_n)^j = \sum_{k=0}^j (-1)^k a_{k,j} (z_nw_n)^k Q_n^{j-k}$, where the 
$a_{k,j}$ are just ordinary binomial coefficients. If we substitute this in 
$\phi$, we obtain
$$
\phi = \sum_{j=0}^{l\wedge m} \sum _{k=0}^j d_ja_{k,j} z_n^{l-j}w_n^{m-j} (z_nw_n)^k Q_n^{j-k}.
$$
Now, after applying Corollary 2.2.7 and changing the summations we obtain $(2.2.5)$. The 
converse statement in the proposition is obvious.\sq
\medskip
Similarly one can prove 
\smallskip
\proclaim\nofrills {\bf Proposition 2.2.10}: Suppose $\phi\in \Z (l,m)$. Then 
$\phi$ is 
${\Cal U}_q(n-2)$-invariant if and only if $\phi$ is of the form
$$
\phi = \sum_{j=0}^{l\wedge m}\sum_{r=0}^{l-j}\sum_{s=0}^{m-j} a_{j,r,s} z_n^{l-j-r}w_n^{m-j-s}
z_{n-1}^rw_{n-1}^s Q_n^j.
$$
\endproclaim
\bigskip

{\bf 3. The quantized algebra of polynomials on the sphere in $\C^n$.}
\medskip
In this chapter we construct the quantized algebra of polynomials on the sphere $S^{2n-1}$ 
from the algebra $\Z$ by putting the invariant element $Q_n$ equal to 1. Furthermore we 
construct an invariant functional on this algebra, give its irreducible decomposition into 
$\U_n$-modules and recover the $\UU$-invariant elements, the so-called zonal spherical 
functions, as $q$-disk polynomials. Finally we obtain an abstract addition formula for these 
$q$-disk polynomials.
\medskip
{\bf 3.1. Definition of ${\widetilde {\Z}}$ and invariant functional.}
\medskip
Consider the following projection map:
$$
\pi: \Z\to \Z /(Q_n -1) =: {\widetilde {\Z}}.
$$
Denote the images of the generators $z_i,w_i$ of $\Z$ under $\pi$ by 
the same symbols, and define $\d$ on those images as in $(2.2.3)$.
This gives a well-defined coaction of $\A_n$ on ${\widetilde {\Z}}$, since $Q_n$ 
is a trivial element for the $\A_n$-coaction. 
So we have the following commutative diagram
$$
\CD
\Z @>\pi>> \widetilde {\Z} \\
@V{\delta}VV    @V{\delta}VV \\
\Z\otimes \A_n @>\pi\otimes id>> \widetilde {\Z}\otimes \A_n.
\endCD \tag 3.1.1
$$
\smallskip
The algebra ${\widetilde {\Z}}$ plays the role of quantized polynomial algebra on the 
$(2n-1)$-sphere $S^{2n-1}$ within ${\C}^n$, and was introduced in [RTF]. It is the same as the 
algebra $A(K\bs G)$ of [NYM, section 4.1].
\smallskip
Recall from section 1.2 that we call a linear functional $h :\widetilde{\Z}\to\C$ $\A_n$-invariant 
if $(h\ot id)\delta(\phi) = h(\phi)1$ for all $\phi\in{\widetilde{\Z}}$. This will imply, see 
$(2.2.4)$, that $h(X.\phi) = \ep(X)h(\phi)$ for $\phi\in{\widetilde {\Z}}$ and $X\in \U_n$. 
In other words, $h$ will be $\U_n$-invariant.
\smallskip
\proclaim\nofrills {\bf Proposition 3.1.1} : On ${\widetilde {\Z}}$ we have a unique normalized,
positive definite $\U_n$-invariant 
functional $h_n :{\widetilde {\Z}} \to {\C}$ given by 
$$
\eqalign{
h_n(z^{\lambda}w^{\mu}) = \delta_{\lambda\mu} &q^{-2((n-1)\lambda_1 + (n-2)\lambda_2 +\hdots 
+\lambda_{n-1})} \times \cr
&{{(q^{-2};q^{-2})_{\lambda_1}\hdots (q^{-2};q^{-2})_{\lambda_n} (q^{-2};q^{-2})_{n-1}}\over 
{(q^{-2};q^{-2})_{|\lambda | +n-1}}}.\cr}\tag 3.1.2
$$
\endproclaim
The proof of this is analogous to [NYM] Proposition 4.5, and uses Proposition 2.2.4.
\smallskip
{\bf Remark}: One can write $(3.1.2)$ equivalently as
$$
\eqalign{
 h_n(z^{\lambda}w^{\mu}) = \delta_{\lambda\mu} 
&q^{|\lambda |^2 + \sum_{i=1}^n (2(i-1)\lambda_i- \lambda_i^2)} \times\cr
&{{(q^2;q^2)_{\lambda_1}\hdots (q^2;q^2)_{\lambda_n} (q^2;q^2)_{n-1}}\over {(q^2;q^2)_{|\lambda 
| +n-1}}}\cr}\tag 3.1.3
$$
when using that $(a^{-1};q^{-1})_m = (-1)^m a^{-m} q^{-{1\over 2} m(m-1)} (a;q)_m$.
\medskip
\proclaim\nofrills {\bf Lemma 3.1.2}: For all $\a,\b\in\z+$ and any continuous function $f$ there holds
$$
\int_0^1 f(q^{-\b}x) x^{\a} (x;q^{-1})_{\b} d_qx = q^{\b(\a+1)} \int_0^1 f(x) x^{\a} (qx;q)_{\b} d_qx.
$$
\endproclaim
\medskip
{\sl Proof}: We have 
$$
\eqalign{
\int_0^1 f(q^{-\b}x) x^{\a} (x;q^{-1})_{\b} d_qx &= (1-q)\sum_{k=0}^{\infty} f(q^{k-\b}) q^{k\a} (q^k;q^{-1})_{\b} q^k\cr
&= (1-q)\sum_{k=\b}^{\infty} f(q^{k-\b}) q^{k\a} (q^k;q^{-1})_{\b} q^k\cr
&= (1-q)\sum_{k=0}^{\infty} f(q^k) q^{(k+\b)\a} (q^{k+\b};q^{-1})_{\b} q^{k+\b}\cr
&= q^{\b(\a+1)} (1-q)\sum_{k=0}^{\infty} f(q^k) q^{k\a} (q^{k+1};q)_{\b} q^k\cr
&= q^{\b(\a+1)} \int_0^1 f(x) x^{\a} (qx;q)_{\b} d_qx.\cr}
$$
\sq
\medskip
Consequently we get
$$
\int_0^1 x^{\a} (x;q^{-1})_{\b} d_q x = q^{\b(\a+1)} (1-q) {{(q;q)_{\a} (q;q)_{\b}}\over 
(q;q)_{\a+\b+1}}
$$
since the right-hand side in Lemma 3.1.2 just represents the $q$-beta 
integral in case $f=1$:
$$
\int_0^1 x^{\a} (qx;q)_{\b} d_qx = (1-q) {{(q;q)_{\a} (q;q)_{\b}}\over 
(q;q)_{\a+\b+1}}.
$$
\medskip
\proclaim\nofrills {\bf Proposition 3.1.3}: In $\wi{\Z}$ we have that span$\{ z^{\lambda}w^{\lambda} 
| \lambda\in {\Bbb Z}^n_{\geq 0} \} = {\C}[Q_1, \hdots, Q_{n-1}]$, and for each 
$\phi=\phi(Q_1,
\hdots,Q_{n-1}) \in {\C}[Q_1, \hdots, Q_{n-1}]$ the value of the Haar integral 
on $\phi$ is 
given by the following multiple Jackson integral:
$$
\eqalign{
h_n(\phi) = &{{(q^2;q^2)_{n-1}}\over {(1-q^2)^{n-1}}} \times \cr
&\int_0^1 \int_0^{Q_{n-1}}\hdots \int_0^{Q_2} \phi(Q_1,\hdots,Q_{n-1}) d_{q^2}Q_1\hdots 
d_{q^2}Q_{n-2} d_{q^2}Q_{n-1}.\cr}\tag 3.1.4
$$
\endproclaim
\smallskip
{\sl Proof}: The first statement follows from Lemma 2.1.2. Hence we only have to verify that $(3.1.4)$ 
is true for any monomial $z^{\lambda}w^{\lambda}$. From Lemma 2.1.2 we get
$$
z^{\lambda}w^{\lambda}= Q_1^{\lambda_1}Q_2^{\lambda_2}\hdots Q_{n-1}^{\lambda_{n-1}}  
({{Q_1}\over {Q_2}};q^{-2})_{\lambda_2}({{Q_2}\over {Q_3}};q^{-2})_{\lambda_3}\hdots 
(Q_{n-1};q^{-2})_{\lambda_n}.
$$
Now substitute this into the right-hand side of $(3.1.4)$, and use that for the general 
$q$-integral one has
$$
\int_0^c f({x\over c})d_qx = c\int_0^1 f(x)d_qx.\tag 3.1.5
$$
Then, by successive use of Lemma 3.1.2, one checks that $(3.1.4)$ yields the same as 
$(3.1.3)$ for $\phi=z^{\lambda}w^{\lambda}$.\sq
\medskip
Finally, we have the following result:
\smallskip
\proclaim\nofrills {\bf Proposition 3.1.4}: The $\ast$-algebra homomorphism 
$\Psi : 
{\widetilde {\Z}}\to \A_n,\ z_k\mapsto t_{nk}$ is well-defined, intertwines $\d$ and $\Delta$
and is injective. So we can apply the results 
of section 1.2 with $A= \A_n$ and $Z= {\widetilde{\Z}}$.
\endproclaim
\smallskip
{\sl Proof}: It is straightforward to verify that $\Psi$ is well-defined and intertwines $\d$
and $\Delta$. Now let us first view $\Psi$ as a map from $\Z$ to 
$A(Mat_q(n))\ot \C [det_q^{-1}]$ (cf. notation in [NYM, section 1.1.]),
where we
assume $det_q^{-1}$ to be a central element but do not assume the identities
$det_q det_q^{-1}=1=det_q^{-1} det_q$ to hold. Suppose $\phi$ is a monomial in 
the elements $t_{ij}\ 
(1\leq i,j\leq n)$ which contains $a_{ij}$ factors $t_{ij}$ $(1\leq i,j\leq n)$. Again we can arrange
things in such a way that $\phi$ has leading term $t^A = t_{11}^{a_{11}}\hdots t_{1n}^{a_{1n}}
t_{21}^{a_{21}}\hdots t_{2n}^{a_{2n}}\hdots t_{nn}^{a_{nn}}$ (so we use the 
following total ordering on the generators $t_{ij}$: $t_{ij}
\prec t_{kl}$ if $i<k$ or $i=k$ and $j<l$ ). We know from [Ko1, Thm. 3.1]  
that monomials in the $t_{ij}$ 
corresponding to different matrices are linearly independent in $A(Mat_q(n))$.
Hence, since we know a linear basis for $\Z$, we must show that the matrix $A(\lambda,\mu)$, 
corresponding to $\Psi (\zw)$, is different for different choices of the pair $(\lambda, \mu)$. 
For this, it suffices to look at the leading term of $\Psi (\zw)$.\hb
Recall that $t_{ij}^{\ast} = S(t_{ji}) = (-q)^{j-i} det_q^{-1} D_{\hat i\hat j}$, 
where $D_{\hat i\hat j}$ is the quantum minor-determinant corresponding to the sub-matrix 
of the matrix $(t_{ij})$ obtained by deleting the $i$-th row and the $j$-th column. The 
leading term of $D_{\hat i\hat j}$ is equal to a constant times $t_{11}\hdots 
t_{j-1,j-1}t_{j,j+1}\hdots t_{i-1, i}t_{i+1,i+1}\hdots t_{nn}$ if $i\geq j$. Hence the 
leading term of $\Psi (\zw) = t_{n1}^{\lambda_1}\hdots t_{nn}^{\lambda_n}(t_{nn}^{\ast})^{\mu_n}\hdots 
(t_{n1}^{\ast})^{\mu_1} = const. det_q^{- |\mu |} t_{n1}^{\lambda_1}\hdots 
t_{nn}^{\lambda_n} 
D_{\hat n\hat n}^{\mu_n}\hdots D_{\hat n\hat 1}^{\mu_1}$ equals a constant times 
$ det_q^{- |\mu |} t^{A(\lambda,\mu)}$, where
$$
\eqalign{ 
A(\lambda,\mu) =&
\pmatrix 0& \hdots &0 \\ \vdots& \ddots &\vdots \\ 0&\hdots &0\\
\lambda_1 &\hdots&\lambda_n \endpmatrix + \cr
 & \pmatrix |\mu |-m_1 &m_1 &0 & &\hdots &0 \\
          0&|\mu | - m_2&m_2&0&\hdots &0\\
          \vdots& &\ddots&\ddots& &\vdots \\
          \vdots& & &\ddots&\ddots&\vdots\\
          0& &\hdots& &|\mu | - m_{n-1}&m_{n-1}\\
          0& &\hdots & & &0\\
\endpmatrix \cr
=&\pmatrix  |\mu |-m_1 &m_1 &0 & &\hdots &0 \\
             0&|\mu | - m_2&m_2&0&\hdots &0\\
             \vdots& &\ddots&\ddots& &\vdots \\
             \vdots& & &\ddots&\ddots&\vdots\\
             0& &\hdots& &|\mu | - m_{n-1}&m_{n-1}\\
             \lambda_1 & &\hdots & & &\lambda_n\\
             \endpmatrix .\cr}
$$
Here we used the short-hand notation $m_k = \mu_1 +\hdots +\mu_k$.
Obviously, different pairs $(\lambda, \mu)$ yield different matrices $A(\lambda,\mu)$, which shows the
linear independence. Finally, using the identity $\sum_{k=1}^n t_{nk}(-q)^{k-n}
D_{\hat n\hat k}=det_q$ in $A(Mat_q(n))$ ([Ko1, (2.10)], [NYM, (1.15.b)]), we find 
$\Psi(Q_n -1) = det_q det_q^{-1}-1$. This shows that $\Psi$ extends to an 
injective homomorphism ${\widetilde{\Z}}\to \A_n$.\sq
\medskip
{\bf Remark} : From this it follows that the algebra $A(K\bs G)$ 
from [NYM, section 4.1] has no additional relations to the relations (4.9.a-d) (loc. cit.). This fact 
was already noted in [VS, Thm. 4.4].
\bigskip

{\bf 3.2 Irreducible decomposition.}
\medskip
With the invariant functional $h_n$ one can define an invariant inner product on ${\widetilde 
{\Z}}$ as follows:
$$
\eqalign{
\la .,.\ra: {\widetilde {\Z}}\times {\widetilde {\Z}} \to \C\cr 
\la\phi,\psi\ra := h_n (\psi^{\ast}\phi).\cr}\tag 3.2.1
$$
This non-degenerate bilinear form then satisfies $\la X.\phi, \psi\ra = \la\phi, X^{\ast}.
\psi\ra$ for all $X\in \U_n$, and all $\phi,\psi\in {\widetilde {\Z}}$. So in particular 
$\la q^h .\phi, \psi\ra = \la \phi, q^h.\psi\ra$ for all $h\in P^{\ast}$ and all $\phi,\psi\in 
{\widetilde {\Z}}$.\hb
Put ${\widetilde {\Z}}(l,m) := \pi (\Z (l,m))$, and let $\HTn(l,m)$ be the orthogonal 
complement with respect to the inner product $(3.2.1)$ of ${\widetilde {\Z}}(l-1,m-1)$ in 
${\widetilde {\Z}}(l,m)$ (recall that we had $Q_n \Z (l-1,m-1)\subset \Z (l,m)$, hence 
${\widetilde {\Z}}(l-1,m-1)\subset {\widetilde {\Z}}(l,m)$). So we have the orthogonal 
direct sum decomposition
$$
{\widetilde {\Z}}(l,m) = {\widetilde {\Z}}(l-1,m-1) \oplus \HTn(l,m).
$$
\smallskip
\proclaim\nofrills {\bf Lemma 3.2.1}: The mapping $\pi : \Z (l,m)\to {\widetilde {\Z}}(l,m)$ 
is injective.
\endproclaim
\medskip
{\sl Proof}: This is a consequence of Proposition 2.1.1 and the fact the the element $Q_n -1$ 
is not homogeneous.\sq
\medskip
From this lemma it follows that $\pi : \Z (l,m)\to {\widetilde {\Z}}(l,m)$ is an isomorphism, 
and hence we obtain from Proposition 2.1.1 that
$$
\text{dim } \wi{\Z} (l,m) = \text{dim } \Z(l,m) = \pmatrix l+n-1\\ n-1 \endpmatrix 
\pmatrix m+n-1\\ n-1 \endpmatrix\tag 3.2.2
$$
and 
$$
\eqalign{
d_n(l,m) := \text{dim } \HTn(l,m) &= \text{dim } {\widetilde {\Z}}(l,m) - \text{dim } {\widetilde 
{\Z}}(l-1,m-1)\cr
&= {{(l+m+n-1)(l+n-2)!(m+n-2)!}\over {l!m!(n-1)!(n-2)!}}.\cr}\tag 3.2.3
$$
\medskip
Moreover we have the decomposition
$$
{\widetilde{\Z}} = \sum_{l,m\in\z+} {\widetilde{\Z}}(l,m).\tag 3.2.4
$$
If we write $\Hn(l,m)$ for the inverse image of $\HTn(l,m)$ under $\pi : \Z (l,m)\to \wi{\Z}(l,m)$, 
there is the direct sum decomposition
$$
\Z(l,m) = Q_n\Z(l-1,m-1) \oplus \Hn(l,m).
$$
\medskip 
\proclaim\nofrills {\bf Proposition 3.2.2}: We have the orthogonal decomposition into inequivalent, 
irreducible ${\Cal U}_q(n)$-modules
$$
{\widetilde {\Z}}(l,m) = \bigoplus_{k=0}^{l\wedge m} \HTn(l-k,m-k).\tag 3.2.5
$$
\endproclaim
\smallskip
{\sl Proof} : It is clear, by the definition of the spaces $\HTn(r,s)$, that ${\widetilde {\Z}}(l,m)$ 
allows the orthogonal direct sum decomposition $(3.2.5)$. Irreducibility of the spaces $\HTn(l-k,m-k)$ 
follows from Proposition 1.2.1:
each non-trivial $\U_n$-invariant subspace of a given $\HTn(l-k,m-k)$ should contain at least one 
$\UU$-invariant element. But there are, according to Proposition 2.2.9, only $l\wedge m +1$ 
linearly independent invariant elements in the space ${\widetilde {\Z}}(l,m)$. Hence none of 
the spaces $\HTn(l-k,m-k)$ 
contains a non-trivial invariant subspace.\hb
 To prove inequivalence, assume that $\HTn(l-k,m-k)
\simeq \HTn(l-k\pr,m-k\pr)$ for $k\neq k\pr$. So in particular $d_n(l-k,m-k)=d_n(l-k\pr,m-k\pr)
=:N$. Take orthonormal bases $\{e_i\}_{i=1}^N$ and $\{f_j\}_{j=1}^N$ in the respective spaces and
construct the elements $\ze_0 =\sum_{i=1}^N \ep(e_i^{\ast})e_i$ and $\ze_0\pr = \sum_{j=1}^N
\ep(f_j^{\ast})f_j$ as in the proof of Proposition 1.2.1. By linear independence it follows that not
all of the $\ep(e_i^{\ast})$ and $\ep(f_j^{\ast})$ can be zero.\hb
Furthermore, put $\eta = \sum_{k=1}^N e_k\ot f_k$. Again we will have that $\Delta\boxtimes\Delta
(\eta) = \eta\ot 1$ and $\Delta(\tilde{\eta})=\tilde{\eta}\ot 1$ where $\tilde{\eta}=\sum_{k=1}^N
e_kf_k$ (cf. Proposition 1.2.1). From this and the orthogonality of the spaces $\HTn(l-k,m-k)$
and $\HTn(l-k\pr,m-k\pr)$ we obtain $\ep(\tilde{\eta}) = h_n(\tilde{\eta})=0$.\hb
Finally, put $\eta_0 = \sum_{k=1}^N \ep(f_k^{\ast})e_k\in \HTn(l-k,m-k)$. Then $\eta_0\neq 0$,
since not all of the $\ep(f_k^{\ast})$ are zero and the $e_i$ are linearly 
independent. Moreover $\ep(\eta_0)=\ep(\tilde{\eta})=0$.
But this means that we have two linearly independent invariant elements within $\HTn(l-k,m-k)$,
namely $\ze_0$ and $\eta_0$ (linearly independent since $\ep(\ze_0)>0$ and $\ep(\eta_0)=0$).
This gives a contradiction, hence the two spaces $\HTn(l-k,m-k)$ and $\HTn(l-k\pr,
m-k\pr)$ cannot be equivalent.\sq
\medskip
{\bf Remark }: Using exactly the same argument as in the proof of the
previous proposition, one shows that the modules $\Hn (l,m)$ are inequivalent
for different choices of the pair $(l,m)$.\hb
From the proof of Proposition 3.2.2 we immediately obtain that each $\HTn(l,m)$ 
contains a unique, up to constants, $\UU$-invariant element. It is called a {\it (zonal) spherical 
function}, or just a spherical element. 
\medskip
\proclaim\nofrills {\bf Lemma 3.2.3}: If $(l,m)\neq (l\pr, m\pr)$, then $\HTn(l,m)\perp 
\HTn(l\pr,m\pr)$.
\endproclaim
\smallskip
{\sl Proof}: Suppose first that $l^{\prime} - l = m^{\prime} - m$, and say this is non-negative. 
This means that $l^{\prime} = l+k,  m^{\prime}= m + k$ for some $k\geq 0$. But then $\HTn(l,m)$ 
and $\HTn(l\pr,m\pr)$ are contained in the same space ${\widetilde{\Z}}(l\pr,m\pr)$, and hence 
orthogonal by Proposition 3.2.2. \hb
On the other hand, if $l^{\prime} - l \neq m^{\prime} - m$, pick $\phi\in \HTn(l,m)$ and 
$\psi\in \HTn(l\pr,m\pr)$. By invariance of the inner product we have
$$
\eqalign{
q^{l-m} \la \phi, \psi\ra & = \la q^{\ep_1 +\hdots+\ep_n}.\phi,\psi\ra\cr
&= \la \phi, q^{\ep_1 +\hdots+\ep_n}.\psi\ra\cr
&= q^{l\pr -m\pr} \la \phi,\psi\ra.\cr}
$$
Since we assumed that $l^{\prime} - l \neq m^{\prime} - m$ and since $q$
is not a root of unity, this gives $\la \phi,\psi\ra = 0$.\sq
\medskip
\proclaim\nofrills {\bf Corollary 3.2.4}: There is the orthogonal, irreducible decomposition 
into inequivalent ${\Cal U}_q(n)$-modules
$$
{\widetilde {\Z}} = \bigoplus_{l,m\in \z+} \HTn(l,m).
$$
\endproclaim
\medskip
{\sl Proof}: This now follows from the decomposition $(3.2.4)$ together with Proposition 3.2.2 
and Lemma 3.2.3 (see also the remark following Proposition 3.2.2).\sq
\bigskip

{\bf 3.3 Zonal spherical functions.}
\medskip
Let us write $\psi (l,m)$ for a spherical element contained in $\HTn(l,m)$, which is unique up 
to constants (see Remark after Proposition 3.2.2). Suppose now that $(l,m)\neq (l\pr, m\pr)$, 
and assume that $l-m = l\pr -m\pr =\b\geq 0$. From Proposition 2.2.9 and Lemma 2.1.2 we know that
$$
\eqalign{
\psi(l,m) = z_n^{l-m}\sum_{j=0}^m c_j z_n^{m-j}w_n^{m-j} = z_n^{\b} p_m(Q_{n-1})\cr
\psi(l\pr,m\pr) = z_n^{l\pr-m\pr}\sum_{j=0}^m c_j z_n^{m\pr-j}w_n^{m\pr-j} = z_n^{\b} p_{m\pr}
(Q_{n-1})\cr}
$$
for certain polynomials $p_m, p_{m\pr}$ of degree $m$ and $m\pr$ respectively. As a consequence 
of the orthogonality of the spaces $\HTn(l,m)$ and $\HTn(l\pr,m\pr)$ we have
$$
\eqalign{
0 &= \la \psi(l\pr,m\pr), \psi(l, m)\ra = h_n(p_m(Q_{n-1})^{\ast} w_n^{\b}z_n^{\b} p_{m\pr}(Q_{n-1}))\cr
&= h_n(p_m(Q_{n-1})^{\ast} p_{m\pr}(Q_{n-1}) (q^2Q_{n-1};q^2)_{\b})\cr}
$$
(cf. Lemma 2.1.2).
Now note that
$$
\int_0^{Q_{n-1}}\hdots \int_0^{Q_2} d_{q^2}Q_1\hdots d_{q^2}Q_{n-2} =
{{(1-q^2)^{n-2}}\over {(q^2;q^2)_{n-2}}} Q_{n-1}^{n-2},\tag 3.3.1
$$
so we obtain from Proposition 3.1.3
$$
0 = \int_0^1 \ov{p_m(Q_{n-1})} p_{m\pr}(Q_{n-1}) Q_{n-1}^{n-2} (q^2Q_{n-1};q^2)_{\b} d_{q^2}Q_{n-1}
\qquad (m\neq m\pr).
$$
But letting $m$ and $m\pr$ vary over $\z+$, these are exactly the orthogonality relations for 
the little $q$-Jacobi polynomials $P^{(n-2,\b)}_m(Q_{n-1};q^2)$ with upper parameters $n-2$ and 
$\b$. In other words, there exist constants $c_m\in {\C}$ such that
$$
p_m(Q_{n-1}) = c^{}_m P^{(n-2,\b)}_m (Q_{n-1};q^2)\qquad\qquad (m\in \z+).
$$
Thus we obtain that in case $l\geq m$ a general spherical element in $\HTn(l,m)$ is given by a 
constant multiple of
$$
z_n^{l-m} P^{(n-2,l-m)}_m (Q_{n-1};q^2)\qquad (l\geq m).
$$
\smallskip
Similarly one can consider the case where $m-l = m\pr-l\pr =\b\geq 0$:
$$
\eqalign{
\psi(l,m) = s_l (Q_{n-1}) w_n^{\b}\cr
\psi(l\pr, m\pr) = s_{l\pr} (Q_{n-1}) w_n^{\b}\cr}
$$
for certain polynomials $s_l,s_{l\pr}$ of degree $l$ and $l\pr$ respectively. In this case the 
orthogonality of the spaces $\HTn(l,m)$ and $\HTn(l\pr,m\pr)$ will yield
$$
0 = h_n(s_l(q^{-2\b}Q_{n-1})^{\ast} s_{l\pr}(q^{-2\b}Q_{n-1}) (Q_{n-1};q^{-2})_{\b})
$$
and thus, using $(3.3.1)$
$$
0=\int_0^1 \ov{s_l(q^{-2\b}Q_{n-1})} s_{l\pr}(q^{-2\b}Q_{n-1}) Q_{n-1}^{n-2} (Q_{n-1};q^{-2})_{\b} 
d_{q^2}Q_{n-1}.
$$
Here ones uses the equalities $z_n Q_{n-1} = q^{-2} Q_{n-1}z_n,\ Q_{n-1}w_n = q^{-2} w_n Q_{n-1}$, 
and $z_n^{\b}w_n^{\b} = (Q_{n-1};q^{-2})_{\b}$ (Lemma 2.1.2).
\smallskip
Observe that Lemma 3.1.2 gives us that 
$$
0 = \int_0^1 \ov{s_l (Q_{n-1})} s_{l\pr}(Q_{n-1}) Q_{n-1}^{n-2} (q^2Q_{n-1};q^2)_{\b} d_{q^2} Q_{n-1}
$$
which gives 
$$
s_l ( Q_{n-1}) = d^{}_l P^{(n-2,\b)}_l (Q_{n-1};q^2)\qquad\quad (l\in {\Bbb N})
$$
for certain $d^{}_l\in \C$.
In other words, the spherical elements of $\HTn(l,m)$ for $l\leq m$ are constant multiples of
$$
P^{(n-2,m-l)}_l (Q_{n-1};q^2) w_n ^{m-l} \qquad (l\leq m).
$$
Summarizing we have
\smallskip
\proclaim\nofrills {\bf Theorem 3.3.1}: For arbitrary $l,m\in \z+$, the $\Cal U_q(n-1)$-invariant 
(i.e. zonal spherical) elements in $\HTn(l,m)$ are constant multiples of
$$
R_{l,m}^{(\a)}(z_n,w_n;q^2) =\left\{\aligned z_n^{l-m} P_m^{(n-2,l-m)}
(Q_{n-1};q^2) \qquad (l\geq m)\\
P_l^{(n-2,m-l)}(Q_{n-1};q^2) w_n^{m-l}\qquad (l\leq m).
\endaligned\right.
$$
\endproclaim
\medskip
{\bf Remark } : This result was already obtained in [NYM, Thm. 4.7].
\medskip
Finally we will calculate the norms of these spherical elements, since we will need them 
later on.\hb
First let us assume that $l-m=\b\geq 0$. Recall that 
$$
\eqalign{
&\int_0^1 R_{l,m}^{(n-2)}(z_n,w_n;q^2)^{\ast} R_{l,m}^{(n-2)}(z_n,w_n;q^2) Q_{n-1}^{n-2}
d_{q^2}Q_{n-1} =\cr
&\int_0^1 P_m^{(n-2,\b)}(Q_{n-1};q^2) P_m^{(n-2,\b)}(Q_{n-1};q^2) Q_{n-1}^{n-2} (q^2Q_{n-1};q^2)_{\b} 
d_{q^2}Q_{n-1} =\cr
&{{(1-q^2)q^{2m(\a+1)}}\over {1-q^{2(\a+\b+2m+1)}}} {{(q^2;q^2)_m (q^2;q^2)_{\b+m}}\over 
{(q^{2(\a+1)};q^2)_m (q^{2(\a+1)};q^2)_{\b+m}}}.\cr}
$$
Thus
$$
\eqalign{
\la &R_{l,m}^{(n-2)} (z_n,w_n;q^2), R_{l,m}^{(n-2)} (z_n,w_n;q^2)\ra \cr
&={{(q^2;q^2)_{n-1}}\over {(1-q^2)^{n-1}}}\times\cr 
&\ \ \int_0^1\int_0^{Q_{n-1}}\hdots\int_0^{Q_2} 
R_{l,m}^{(n-2)} (z_n,w_n;q^2)^{\ast} R_{l,m}^{(n-2)} (z_n,w_n;q^2)d_{q^2}Q_1
\hdots d_{q^2}Q_{n-1}\cr
&={{(q^2;q^2)_{n-1}}\over {(1-q^2)^{n-1}}} {{(1-q^2)^{n-2}}\over 
{(q^2;q^2)_{n-2}}}\times\cr
&\ \ \int_0^1 P_m^{(n-2,\b)}(Q_{n-1};q^2) P_m^{(n-2,\b)}(Q_{n-1};q^2) Q_{n-1}^{n-2} 
(q^2Q_{n-1};q^2)_{\b} d_{q^2}Q_{n-1}\cr
&={{1-q^{2(n-1)}}\over {1-q^2}}\times\cr
&\ \ \int_0^1 P_m^{(n-2,\b)}(Q_{n-1};q^2) P_m^{(n-2,\b)}(Q_{n-1};q^2) Q_{n-1}^{n-2} 
(q^2Q_{n-1};q^2)_{\b} d_{q^2}Q_{n-1}\cr
&={{(1-q^{2(n-1)})q^{2m(n-1)}}\over {1-q^{2(n+l+m-1)}}}  {{(q^2;q^2)_l (q^2;q^2)_{m}}\over 
{(q^{2(n-1)};q^2)_l (q^{2(n-1)};q^2)_{m}}}.\cr}
$$
Similarly, in case $m-l=\b\geq 0$ we saw that (cf. Lemma 3.1.2)
$$
\eqalign{
\la &R_{l,m}^{(n-2)} (z_n,w_n;q^2), R_{l,m}^{(n-2)} (z_n,w_n;q^2)\ra \cr
&={{(q^2;q^2)_{n-1}}\over {(1-q^2)^{n-1}}} {{(1-q^2)^{n-2}}\over 
{(q^2;q^2)_{n-2}}}\times\cr
&\ \ \int_0^1 P_l^{(n-2,\b)}(q^{-2\b}Q_{n-1};q^2) P_l^{(n-2,\b)}(q^{-2\b}Q_{n-1};q^2) Q_{n-1}^{n-2} 
(Q_{n-1};q^{-2})_{\b} d_{q^2}Q_{n-1}\cr
&={{1-q^{2(n-1)}}\over {1-q^2}}{q^{2\b(n-1)}}\times\cr
&\ \ \int_0^1 P_l^{(n-2,\b)}(Q_{n-1};q^2) P_l^{(n-2,\b)}(Q_{n-1};q^2) Q_{n-1}^{n-2} 
(q^2Q_{n-1};q^2)_{\b} d_{q^2}Q_{n-1}.\cr}
$$
Thus, in case $m-l=\b\geq 0$ we get
$$
\eqalign{
\la R_{l,m}^{(n-2)} (z_n,w_n;q^2), R_{l,m}^{(n-2)} (z_n,w_n;q^2)\ra = &
{{(1-q^{2(n-1)})q^{2l(n-1)}q^{2\b(n-1)}}\over {1-q^{2(n+l+m-1)}}} \times\cr
&{{(q^2;q^2)_l (q^2;q^2)_{m}}\over {(q^{2(n-1)};q^2)_l (q^{2(n-1)};q^2)_{m}}}\cr
=&{{(1-q^{2(n-1)})q^{2m(n-1)}}\over {1-q^{2(n+l+m-1)}}}\times\cr &{{(q^2;q^2)_l 
(q^2;q^2)_{m}}\over {(q^{2(n-1)};q^2)_l (q^{2(n-1)};q^2)_{m}}}.\cr}
$$
So we obtain
\smallskip
\proclaim\nofrills {\bf Proposition 3.3.2} : If $l,m\in\z+$ are arbitrary, then the norm 
of the polynomial $R_{l,m}^{(\a)} (z_n,w_n;q^2)$ is given by
$$
\eqalign{
\Vert R_{l,m}^{(\a)} (z_n,w_n;q^2)\Vert &= h_n\bigl(R_{l,m}^{(\a)} (z_n,w_n;q^2)^{\ast} 
R_{l,m}^{(\a)} (z_n,w_n;q^2)\bigr)^{1\over 2}\cr
&= (c_{l,m}^{(\a)})^{1\over 2}\cr}
$$
in which
$$
c_{l,m}^{(\a)} = {{(1-q^{2(\a+1)})q^{2m(\a+1)}}\over {1-q^{2(\a+l+m+1)}}} {{(q^2;q^2)_l 
(q^2;q^2)_{m}}\over {(q^{2(\a+1)};q^2)_l (q^{2(\a+1)};q^2)_{m}}}\tag 3.3.2
$$
and $\a=n-2\in \z+$.
\endproclaim
{\bf Remark }: Note that $c_{l,m}^{(\a)}$ is not symmetric in $l$ and $m$, unlike
in the classical case. We remark that
$$
\eqalign{
h_n(R_{l,m}^{(\a)}(z_n,w_n;q^2) R_{l,m}^{(\a)}(z_n,w_n;q^2)^{\ast}) &= c_{m,l}^{(\a)}\cr
&= h_n(R_{m,l}^{(\a)}(z_n,w_n;q^2)^{\ast} R_{m,l}^{(\a)}(z_n,w_n;q^2))\cr}
$$
since $R_{l,m}^{(\a)}(z_n,w_n;q^2)^{\ast}=R_{m,l}^{(\a)}(z_n,w_n;q^2)$.
\bigskip
{\bf 3.4 Associated spherical functions}.
\medskip
Suppose we are given the $\ast$-algebra ${\Cal Z}_{n-1}$ with generators $z_i\pr, w_i\pr$ 
$(1\leq i\leq n-1)$ and with corresponding projection $\pi\pr : {\Cal Z}_{n-1}\to 
{\widetilde{\Cal Z}}_{n-1}$. Recall that we had the embedding $\iota^{(n-1,n)} : 
{\Cal Z}_{n-1}(r,s)\hookrightarrow \Z(r,s)$ $(r,s\in\z+)$. Then, using the map 
$\pi\circ\iota^{(n-1,n)}\circ (\pi\pr)^{-1}$, one can identify $\phi=\phi(z\pr, w\pr)\in 
{\widetilde{\Cal Z}}_{n-1}(r,s)$ with $Q_{n-1}^{{1\over 2}(r+s)} \phi(z Q_{n-1}^{-{1\over 2}}, 
w Q_{n-1}^{-{1\over 2}})\in {\widetilde {\Z}}(r,s)$.
\smallskip
Given $l,m\in \z+$ and $0\leq r\leq l,\ 0\leq s\leq m$, define the following elements in $\Z(l,m)$:
$$
\eqalign{
\psi(l,m;r,s) 
= &Q_n^{{1\over 2}(l-r+m-s)} R_{l-r,m-s}^{(n-2+r+s)}(z_n Q_n^{-{1\over 2}}, w_n 
Q_n^{-{1\over 2}};q^2)\times\cr &Q_{n-1}^{{1\over 2}(r+s)} R_{r,s}^{(n-3)}(z_{n-1} 
Q_{n-1}^{-{1\over 2}},w_{n-1} Q_{n-1}^{-{1\over 2}};q^2).\cr} 
$$
Denote their restrictions in $\widetilde {\Z}(l,m)$ via $\pi$ by the same symbols.
\smallskip
\proclaim\nofrills {\bf Proposition 3.4.1}: In $\widetilde{\Z}$ there holds $\la 
\psi(l\pr,m\pr;r\pr,s\pr), \psi (l,m ;r,s)\ra = 0$ whenever\hb $(l\pr,m\pr,r\pr,s\pr)\neq (l,m,r,s)$.
Moreover
$$
\la \psi(l,m;r,s), \psi (l,m ;r,s)\ra = {{1-q^{2(\a+1)}}\over {1-q^{2(\a+r+s+1)}}} 
c_{l-r,m-s}^{(\a+r+s)}  c_{r,s}^{(\a-1)}
$$
in which $\a=n-2$.
\endproclaim
\smallskip
{\sl Proof }: This is done by direct calculation. From $(3.1.3)$ we see immediately that 
the above inner product is zero if we do not have $r-s=r\pr-s\pr$ and $l-m=l\pr-m\pr$. 
If we do, then it follows that $R_{l-r,m-s}^{(n-2+r+s)}(z_n,w_n;q^2)^{\ast} 
R_{l\pr-r\pr,m\pr -s\pr}^{(n-2+r\pr+s\pr)}$ is a polynomial in $Q_{n-1}$, and hence (see Lemma 2.1.2)
$$
\eqalign{
&\psi(l,m;r,s)^{\ast}\psi(l\pr,m\pr;r\pr,s\pr) =\cr
&R_{l-r,m-s}^{(n-2+r+s)}(z_n,w_n;q^2)^{\ast}
R_{l\pr-r\pr,m\pr-s\pr}^{(n-2+r+s)}(z_n,w_n;q^2)\times\cr
&Q_{n-1}^{{1\over 2}(r+s+r\pr+s\pr)}
R_{r,s}^{(n-3)}(z_{n-1}\pr,w_{n-1}\pr;q^2)^{\ast} 
R_{r\pr,s\pr}^{(n-3)}(z_{n-1}\pr,w_{n-1}\pr;q^2).\cr}
$$
Here we used the short-hand notation $z_{n-1}\pr=z_{n-1}Q_{n-1}^{-{1\over 2}},
w_{n-1}\pr =w_{n-1}Q_{n-1}^{-{1\over 2}}$. Observing that the third line is a polynomial 
in ${Q_{n-2}}\over {Q_{n-1}}$, since we assumed $r-s=r\pr-s\pr$, we can now calculate 
the value of $h_n$ using $(3.1.4)$, $(3.3.1)$ and $(3.1.5)$:
$$
\eqalign{
&h_n(\psi(l,m;r,s)^{\ast}\psi(l\pr,m\pr;r\pr,s\pr)) =
{{(q^2;q^2)_{n-1}}\over {(1-q^2)^{n-1}}} {{(1-q^2)^{n-3}}\over {(q^2;q^2)_{n-3}}}\times\cr
&\bigl(\int_0^1 R_{r,s}^{(n-3)}(z_{n-1}\pr,w_{n-1}\pr;q^2)^{\ast} 
R_{r\pr,s\pr}^{(n-3)}(z_{n-1}\pr,w_{n-1}\pr;q^2) (Q_{n-2}\pr)^{n-3}d_{q^2}Q_{n-2}\pr\bigr) \times\cr
&\int_0^1  R_{l-r,m-s}^{(n-2+r+s)}(z_n,w_n;q^2)^{\ast}
R_{l\pr-r\pr,m\pr-s\pr}^{(n-2+r+s)}(z_n,w_n;q^2)
Q_{n-1}^{{1\over 2}(r+s+r\pr+s\pr)} Q_{n-1}^{n-2} d_{q^2}Q_{n-1}\cr}
$$
where in the first integral the primed elements are now considered as elements of
$\widetilde {\Cal Z}_{n-1}$. Thus (cf. Proposition 3.3.2) this becomes
$$
\eqalign{
& =\d_{rr\pr}\d_{ss\pr} {{(1-q^{2(n-2)})(1-q^{2(n-1)})}\over {(1-q^2)^2}} {{1-q^2}\over 
{1-q^{2(n-2)}}} c_{r,s}^{(n-3)}\times\cr
&\ \ \int_0^1  R_{l-r,m-s}^{(n-2+r+s)}(z_n,w_n;q^2)^{\ast}
R_{l\pr-r\pr,m\pr-s\pr}^{(n-2+r+s)}(z_n,w_n;q^2)
Q_{n-1}^{n-2+r+s}d_{q^2}Q_{n-1}\cr
& =\d_{ll\pr}\d_{mm\pr}\d_{rr\pr}\d_{ss\pr} {{1-q^{2(n-1)}}\over {1-q^2}} c_{r,s}^{(n-3)} 
{{1-q^2}\over {1-q^{2(n-1+r+s)}}} c_{l-r,m-s}^{(n-2+r+s)}\cr
&=\d_{ll\pr}\d_{mm\pr}\d_{rr\pr}\d_{ss\pr}{{1-q^{2(n-1)}}
\over {1-q^{2(n-1+r+s)}}}c_{l-r,m-s}^{(n-2+r+s)}c_{r,s}^{(n-3)}\cr}
$$
which proves the proposition.\sq
\medskip
\proclaim\nofrills {\bf Proposition 3.4.2}: For fixed $l$ and $m$, $\psi(l,m;r,s)\in 
\HTn(l,m)$ for all $0\leq r\leq l$ and all $0\leq s\leq m$. Moreover, $F\in \HTn(l,m)$ 
is ${\Cal U}_q (n-2)$-invariant if and only if $F\in \text{span} \{\psi(l,m;r,s)\ \vert\ 
0\leq r\leq l,\ 0\leq s\leq m\}$.
\endproclaim
\smallskip
{\sl Proof }: For $0\leq j\leq l\wedge m,\ 0\leq r\leq l-j,\ 0\leq s\leq m-j$ consider the 
elements $Q_n^j \psi (l-j,m-j;r,s)$ in $\Z (l,m)$. Their restrictions to $\tilde {\Z}$ are 
mutually orthogonal by the previous proposition, hence they are linearly independent in $\Z 
(l,m)$. Since they are all ${\Cal U}_q(n-2)$-invariant, they will span the entire space of 
${\Cal U}_q(n-2)$-invariant elements within $\Z(l,m)$ because of their number (cf. 
Proposition 2.2.10). In the same way all elements $Q_n^j \psi (l-j,m-j;r,s)$ with 
$1\leq j\leq l\wedge m$ and  $0\leq r\leq l-j,\ 0\leq s\leq m-j$ will span the subspace 
of ${\Cal U}_q (n-2)$-invariant elements in $Q_n \Z(l-1,m-1)$. Now, using the orthogonality 
of the $\psi(l,m;r,s)$ together with Lemma 1.2.2,
we conclude that the $\psi(l,m;r,s)$ with $0\leq r\leq l,\ 0\leq s\leq m$ are orthogonal 
to the whole of $\tilde{\Z}(l-1,m-1)$, hence belong to $\HTn(l,m)$. Because of their number 
the second part of the proposition is also clear.\sq
\medskip
For given $l,m\in \Bbb N$ and $0\leq r\leq l$, $0\leq s\leq m$ put
$$
\Hn(l,m;r,s) = Q_n^{{1\over 2}(l-r-m+s)} R_{l-r,m-s}^{(n-2+r+s)} (z_n Q_n^{-{1\over 2}}, 
w_n Q_n^{-{1\over 2}};q^2) \Hnn(r,s).
$$
Then clearly $\Hn(l,m;r,s)\subset \Z(l,m)$ ($0\leq r\leq l$, $0\leq s\leq m$).
\smallskip
\proclaim\nofrills {\bf Lemma 3.4.3}: For all $0\leq r\leq l$ and all $0\leq s\leq m$ we 
have the inclusion
$$
\Hn(l,m;r,s)\subset \Hn(l,m).
$$
\endproclaim
\smallskip
{\sl Proof }: As $\UU$-modules there is the isomorphism $\Hn(l,m;r,s)\cong \Hnn(r,s)$. 
Since $\Hnn(r,s)$ is an irreducible $\UU$-module (for it is isomorphic to the module 
$\HTnn(r,s)$), and since  $\psi(l,m;r,s)\in \Hn(l,m;r,s)$, we know that $\UU . \psi(l,m;r,s) 
= \Hn(l,m;r,s)$ for all $0\leq r\leq l$ and all $0\leq s\leq m$. On the other hand, 
$\psi(l,m;r,s)\in \Hn(l,m)$ (Proposition 3.4.2), hence $\UU . \psi(l,m;r,s)\subset \Hn(l,m)$. 
Thus we see that $\Hn(l,m;r,s)\subset \Hn(l,m)$ for all $0\leq r\leq l$ and all $0\leq s\leq m$.\sq
\medskip
\proclaim\nofrills {\bf Proposition 3.4.4}: We have the following direct sum decomposition 
into irreducible ${\Cal U}_q (n-1)$-modules
$$
\Hn(l,m) = \bigoplus_{r=0}^l \bigoplus_{s=0}^m \Hn(l,m;r,s).\tag 3.4.1
$$
\endproclaim
\smallskip
{\sl Proof }: From the previous lemma we get that the direct sum on the right hand side is 
contained in $\Hn(l,m)$. Counting dimensions gives the equality.\sq
\medskip
Write $\HTn(l,m;r,s) = \pi (\Hn(l,m;r,s))$.
\smallskip
\proclaim\nofrills {\bf Proposition 3.4.5}: We have the following orthogonal decomposition 
into irreducible ${\Cal U}_q (n-1)$-modules
$$
\HTn(l,m) = \bigoplus_{r=0}^l \bigoplus_{s=0}^m \HTn(l,m;r,s).\tag 3.4.2
$$
\endproclaim
\smallskip
This proposition follows from the previous one and the following
\smallskip
\proclaim\nofrills {\bf Proposition 3.4.6}: If $\{ \pi\pr (g_i(r,s))\ :\ i=1,\hdots, 
d_{n-1}(r,s)\}$ forms an orthonormal basis for $\HTnn(r,s)$ with respect to the inner 
product defined by $h_{n-1}$, then the set 
$\{( {{1-q^{2(n-1)}}\over {1-q^{2(n+r+s-1)}}}c_{l-r,m-s}^{(n+r+s-2)})^{-{1\over 2}}  
R_{l-r,m-s}^{(n-2+r+s)}(z_n, w_n;q^2) \pi(g_i(r,s))\}$ with 
$0\leq r\leq l$, $0\leq s\leq m$ and $i=1,\hdots,d_{n-1}(r,s)$ forms an orthonormal 
basis for $\HTn(l,m)$ with respect to the inner product on $\tilde{\Z}$ defined by $h_n$. 
\endproclaim
\smallskip
{\sl Proof }: The proof of this is along the same lines as the proof of Proposition 
3.4.1.\sq
\medskip
{\bf Remark }: The elements of Proposition 3.4.6 are called {\it associated spherical elements} 
in $\Hn (l,m)$.
\bigskip

{\bf 3.5 Addition formula for $q$-disk polynomials.}
\medskip
We are now at the stage where we can state the addition formula for the $q$-disk polynomials. 
For this we use the concrete realization of $\widetilde{\Z}$ as a $\ast$-subalgebra of $\A_n$ 
as in section 3.1.\hb
So let us identify $z_i=t_{ni},\ w_i=t_{ni}^{\ast}$. It should be noted that under this 
correspondence the coaction $\d$ is merely the comultiplication $\Delta$. Write $\tau$ 
for the anti-linear involutive algebra automorphism $\ast\circ S$ from $\A_n$ to $\A_n$. 
We know that $\tau(t_{ij}) =t_{ji}$, since $t_{ij}^{\ast} = S(t_{ji}) = 
(-q)^{j-i}\xi_{\hat i\hat j}det_q^{-1}$ where
$\xi_{IJ}$ denotes the quantum minor determinant corresponding to the two subsets
$I,J\subset \{1,\hdots,n\}$ (with $\# (I) =\#(J)$). So in particular we get 
$\tau(z_n) = z_n$. Moreover, recall from [NYM, (3.2)] that
$$
(D_{IJ})^{\ast} = S(D_{JI}) = {{\text{sgn}_q (J;J^c)}\over {\text{sgn}_q(I;I^c)}} 
D_{I^cJ^c} det_q^{-1}\tag 3.5.1
$$
in which $I^c$ denotes the complement of $I$ in $\{1,\hdots ,n\}$ and
$$
\text{sgn}_q (I;J) = \left\{\aligned 0 \qquad I\cap J\neq\emptyset\\
(-q)^{l(I;J)}\qquad I\cap J =\emptyset\endaligned\right.
$$
where $l(I;J) = \# \{(i,j)\in I\times J\ | \ i>j\}$.
Using $(3.5.1)$ we see that $S(t_{nn}^{\ast}) = t_{nn}$, whence $\tau(w_n) = w_n$. So we conclude 
that $\tau(Q_{n-1})=Q_{n-1}$.\hb
As observed in a remark in section 1.2 we can exhibit $\HTn(l,m)$ as the row space 
$A_{\pi(1)}$ for certain $\pi\in\Sigma$ such that $\pi_{11}= R_{l,m}^{(n-2)}(z_n,w_n;q^2)$ (since
it is easily seen that $\ep(R_{l,m}^{(n-2)}(z_n,w_n;q^2))=1$). The basis elements $\{\pi_{1i}\}$ then 
correspond 
to the elements as given in Proposition 3.4.6 with $r+s\neq 0$.
By virtue of that same remark and the fact that for $\pi\in\Sigma$ one has
$(id\ot\tau)\Delta(\pi_{11})=\sum_k \pi_{1k}\ot\pi_{1k}$, we can write
$$
\eqalign{
(id\ot\tau)\Delta R_{l,m}^{(n-2)}(z_n,w_n;q^2) = &h_n(R_{l,m}^{(n-2)}(z_n,w_n;q^2)^{\ast} 
R_{l,m}^{(n-2)}(z_n,w_n;q^2))\times\cr
\sum_{r=0}^l\sum_{s=0}^m \sum_{i=1}^{d_{n-1}(r,s)} (&a_{l,m,r,s}^{-{1\over 2}}
R_{l-r,m-s}^{(n-2+r+s)}(z_n,w_n;q^2) \pi(g_i(r,s)))\ot \cr
&a_{l,m,r,s}^{-{1\over 2}} R_{l-r,m-s}^{(n-2+r+s)}(z_n,w_n;q^2) \pi(g_i(r,s))\cr}
$$
where $a_{l,m,r,s} = {{1-q^{2(n-1)}}\over {1-q^{2(n+r+s-1)}}}c_{l-r,m-s}^{(n+r+s-2)}$. 
Here we choose the orthonormal bases $\{g_i(r,s)\}$ for all $r,s\in \z+$ in such a way 
that $g_1(r,s)= 
(c_{r,s}^{(n-3)})^{-{1\over 2}} Q_{n-1}^{{1\over 2}(r+s)}\times\hb
R_{r,s}^{(n-3)}(z_{n-1} 
Q_{n-1}^{-{1\over 2}} ,w_{n-1}Q_{n-1}^{-{1\over 2}})$. Let us 
pull this equality in $\HTn(l,m)\ot\HTn(l,m)$ back to $\Hn(l,m)\ot \Hn(l,m)$. 
We get
$$
\eqalign{
(Q_n\ot Q_n)^{{1\over 2}(l+m)} R_{l,m}^{(n-2)}&({{(id\ot\tau)\Delta(z_n)}\over {(Q_n\ot Q_n)^{1\over 2}}},
{{(id\ot\tau)\Delta(w_n)}\over {(Q_n\ot Q_n)^{1\over 2}}};q^2) =
 c_{l,m}^{(n-2)}\times\cr
\sum_{r=0}^l\sum_{s=0}^m \sum_{i=1}^{d_{n-1}(r,s)}
 a_{l,m,r,s}^{-1}
&Q_n^{{1\over 2}(l-r+m-s)}
R_{l-r,m-s}^{(n-2+r+s)}(z_nQ_n^{-{1\over 2}},w_nQ_n^{-{1\over 2}};q^2)
g_i(r,s)) \ot\cr 
&Q_n^{{1\over 2}(l-r+m-s)}
R_{l-r,m-s}^{(n-2+r+s)}(z_nQ_n^{-{1\over 2}},w_nQ_n^{-{1\over 2}};q^2) 
g_i(r,s).\cr}\tag 3.5.2
$$
Now recall the projection $\rho^{(n,2)} : \Z \to {\Cal Z}_2$ which puts 
$z_1,\hdots, z_{n-2}, w_1,\hdots, w_{n-2}$ equal to zero. 
Let us write $\g, \d, -q^{-1}\b,\a$ for the generators $z_1\pr,z_2\pr, w_1\pr,
w_2\pr$ of ${\Cal Z}_2$ respectively.
We also will write $D$ for 
the generator of the center of ${\Cal Z}_2$: $D= \a^{\ast}\a + \g\g^{\ast} =\d\a 
-q^{-1}\b\g$. Another way of writing $D$ is $D=\a\a^{\ast} + q^2 \g\g^{\ast} = \a\d - 
q\b\g$. Having this, we can write the projection $\rho^{(n,2)}$ as
$$
\eqalign{
\rho^{(n,2)} : \Z&\to {\Cal Z}_2\cr
z_i,w_i&\to 0\qquad (i=1,\hdots,n-2)\cr
z_{n-1}&\to\g\cr
z_n&\to\d\cr
w_{n-1}&\to\g^{\ast}=-q^{-1}\b\cr
w_n&\to\d^{\ast}=\a.\cr}\tag 3.5.3
$$
\smallskip

\proclaim\nofrills {\bf Lemma 3.5.1} : For $0\leq r\leq l,\ 0\leq s\leq m$ pick a 
basis $\{g_i(r,s)\}$ of $\Hnn(r,s)$ such that in each case $g_1(r,s)= 
(c_{r,s}^{(n-3)})^{-{1\over 2}} Q_{n-1}^{{1\over 2}(r+s)}R_{r,s}^{(n-3)}(z_{n-1}
Q_{n-1}^{-{1\over 2}} ,w_{n-1} 
Q_{n-1}^{-{1\over 2}})$. Then:
$$
\rho^{(n,2)}(g_i(r,s)) = \d_{i1} (c_{r,s}^{(n-3)})^{-{1\over 2}} 
\g^r (\g^{\ast})^s.
$$
\endproclaim
{\sl Proof } : From Proposition 3.4.4 we have the decomposition
$$
\Hnn (r,s) =\bigoplus_{u=0}^r\bigoplus_{v=0}^s \Hnn (r,s;u,v).
$$
This immediately yields that $^{\rho^{(n,2)}}\vert_{\Hnn(r,s;u,v)}\neq 0$ if
and only if  
$(u,v)=(0,0)$. Since $\Hnn(r,s;0,0)$ is the one-dimensional space spanned by 
$g_1(r,s)$, the lemma now follows from an easy computation.\sq
\smallskip
\proclaim\nofrills {\bf Lemma 3.5.2} : $(id\ot\tau)\Delta(z_n)=
\sum_{k=1}^nz_k\ot z_k$\hb
\qquad\qquad $(id\ot\tau)\Delta(w_n)=\sum_{k=1}^n q^{2(n-k)} w_k\ot w_k.$
\endproclaim

{\sl Proof } : The first equality follows directly from $(2.2.3)$ and the fact that 
$\tau(t_{ij})= t_{ji}$. As for second one, use $(3.5.1)$ to obtain $\tau(t_{kn}^{\ast})= 
q^{2(n-k)}w_k$. 
Combining this with $(2.2.3)$ yields the stated result.\sq
\smallskip
\proclaim\nofrills {\bf Corollary 3.5.3} : $(id\ot\rho^{(n,2)})(id\ot\tau)\Delta(z_n) =
z_{n-1}\ot\g + z_n\ot\d$\hb
\hskip4cm $(id\ot\rho^{(n,2)})(id\ot\tau)\Delta(w_n)=q^2w_{n-1}\ot\g^{\ast} + w_n\ot\d^{\ast} 
= -qw_{n-1}\ot\b + w_n\ot\a$.
\endproclaim
{\sl Proof } : Immediate from $(3.5.3)$ and the previous lemma.\sq
\medskip
Consider the $\ast$-algebras ${\Cal X}$ and
${\Cal Y}$ generated by the elements
$$
\alignat 5
& \Cal X :&\quad &Q=Q_n & &\qquad &\qquad &\Cal Y :&\quad 
&D= D\\
& & &X_1 =z_{n-1} & &\qquad & & & &Y_1=\g\\
& & &X_1^{\ast}=w_{n-1} & &\qquad & & & 
&Y_1^{\ast} = \g^{\ast}\tag 3.5.4\\
& & &X_2 =z_n & &\qquad & & & &Y_2 = \d\\
& & &X_2^{\ast}= w_n & &\qquad & & & 
&Y_2^{\ast}=\d^{\ast}
\endalignat
$$
and with $\ast$-structures
$$
\alignat 3
&Q^{\ast} = Q & &\qquad & &D^{\ast}=D\\
&(X_1)^{\ast} = X_1^{\ast} & &\qquad & 
&(Y_1)^{\ast}=Y_1^{\ast}\tag 3.5.5\\
&(X_2)^{\ast}=X_2^{\ast} & &\qquad & 
&(Y_2)^{\ast}=Y_2^{\ast}.
\endalignat
$$
(so we merely changed notations). It is straightforward from $(2.1.1-4)$ 
that the following relations are satisfied:
\proclaim\nofrills {\bf Lemma 3.5.4} : One has:
$$
\eqalign{
X_1X_2 = &q X_2 X_1\cr
X_1^{\ast} X_2 =& qX_2 X_1^{\ast}\cr
X_2^{\ast}X_2= &q^2X_2X_2^{\ast} + 
(1-q^2)Q\cr
X_1^{\ast}X_1 = &q^2 X_1X_1^{\ast}  
+(1-q^2) (Q-X_2X_2^{\ast})\cr
Q\ &\text{central}\cr
Y_1Y_2 = &qY_2Y_1\cr
Y_1^{\ast}Y_2 = &qY_2Y_1^{\ast}\cr
Y_1Y_1^{\ast} = &Y_1^{\ast} Y_1\cr
D =&Y_1Y_1^{\ast} + Y_2Y_2^{\ast}
= q^2Y_1^{\ast}Y_1 + Y_2^{\ast}Y_2.\cr}
\tag 3.5.6
$$
\endproclaim
\smallskip
{\bf Remark } : Centrality of $Q$ in $\Cal X$ does not follow automatically 
from the first four relations above, but is imposed on the algebra $\Cal X$. 
However, $D$ is clearly central.

\medskip
Write $\Cal B={\Cal X}\ot {\Cal Y}$ and identify $X_1$ with $X_1\ot 1
\in \Cal B$, and so on. We will prove that the relations in Lemma 3.5.4 
are in fact the only non-trivial relations in $\Cal B$.\hb
\smallskip
\proclaim\nofrills {\bf Lemma 3.5.5} : A linear basis for $\Cal X$ is given by the set of all 
monomials\hb
 $\{X_1^rX_2^s(X_2^{\ast})^t(X_1^{\ast})^u(Q\pr)^v\ \vert\ r,s,t,u,v\in \z+\}$
where $Q\pr = Q- X_1X_1^{\ast} -X_2X_2^{\ast}$.
\endproclaim
{\sl Proof } : Rewrite the relations in $\Cal X$ in terms of $X_1, X_2,
X_1^{\ast}, X_2^{\ast}$
and $Q\pr$ as follows:
$$
\eqalign{
X_1X_2 &= q X_2 X_1\cr
X_1^{\ast} X_2 &= qX_2 X_1^{\ast}\cr
X_2^{\ast} X_2&= X_2 X_2^{\ast} + (1-q^2)(Q\pr + X_1 X_1^{\ast})\cr
X_1^{\ast} X_1&= X_1 X_1^{\ast} + (1-q^2)Q\pr.\cr}
$$
From this it readily follows that $\Cal X$ is spanned by all the monomials
as in the lemma. Hence we only need to show linear independence of this set.
For this it suffices to show linear independence of the highest order terms,
which, as elements of $\Z$, equal
$$
\eqalign{
\ho (X_1^rX_2^s(X_2^{\ast})^t(X_1^{\ast})^u(Q\pr)^v) &= 
\ho (z_{n-1}^r z_n^s w_n^t w_{n-1}^u (Q_n\pr)^v)\cr
&= z_{n-1}^r z_n^s w_n^t w_{n-1}^u z_{n-2}^v w_{n-2}^v.\cr}
$$
By virtue of Proposition 2.1.1 these are linearly independent as elements
of $\Z$, hence also as elements of $\Cal X$.\sq
\smallskip
From this and Proposition 2.1.1 we get as a consequence
\smallskip
\proclaim\nofrills {\bf Proposition 3.5.6} : A linear basis for $\Cal B$ is
given 
by the set
of all monomials $(X_1\ot 1)^r(X_2\ot 1)^s 
(X_2^{\ast}\ot 1)^t(X_1^{\ast}\ot 1)^u(Q\pr\ot 1)^v 
(1\ot Y_1)^k(1\ot Y_1^{\ast})^l(1\ot Y_2)^m(1\ot Y_2^{\ast})^p$ 
where $k,l,m,p,r,s,t,u,v\in \z+$.
\endproclaim
\smallskip
Now we can prove:
\proclaim\nofrills {\bf Proposition 3.5.7} : The relations $(3.5.6)$ are the only non-trivial 
relations among the generators of $\Cal B$.
\endproclaim
{\sl Proof } : Write $\Cal E$ for the 
$\ast$-algebra 
with abstract generators $X_1,X_1^{\ast},
X_2,X_2^{\ast},Q$ and $Y_1,Y_1^{\ast}, Y_2,Y_2^{\ast},D$ and relations $(3.5.6)$.
Furthermore impose that all of the first five generators commute with all of
the last five ones.
The 
$\ast$-structure on $\Cal E$ is given by $(3.5.5)$. From $(3.5.6)$ we see that $\Cal E$ is 
spanned by all the monomials of the form:\hb $X_1^rX_2^s(X_2^{\ast})^t
(X_1^{\ast})^u (Q\pr)^v
Y_1^k(Y_1^{\ast})^lY_2^m(Y_2^{\ast})^p$ with $k,l,m,p,r,s,t,u,v\in \z+$. 
As before we put $Q\pr = Q- X_1X_1^{\ast} -X_2X_2^{\ast}$.
There is a unique 
surjective $\ast$-algebra homomorphism $\Theta : \Cal E\to \Cal B$ sending $X_1,X_2,X_2^{\ast},
X_1^{\ast},Q\pr,\hb Y_1,Y_1^{\ast},Y_2$ and $Y_2^{\ast}$ to $z_{n-1}\ot 1,z_n\ot 1,w_n\ot 
1,w_{n-1}\ot 1, Q_n\pr\ot 1, 1\ot\g, 1\ot\g^{\ast},1\ot\d$ and $1\ot\d^{\ast}$ respectively. 
It now easily follows from Proposition 3.5.6 that this is actually an isomorphism.\sq
\medskip
With the notation as in $(3.5.4)$ and with the aid of Lemma 3.5.1 and 
Corollary 3.5.3, we 
can write down the effect of applying $id\ot\rho^{(n,2)}$ to $(3.5.2)$. 
It reads:
$$
\eqalign{
R_{l,m}^{(n-2)}(X_1\ot &Y_1 +X_2\ot Y_2, q^2X_1^{\ast}\ot Y_1^{\ast}
+X_2^{\ast}\ot Y_2^{\ast}, QD; q^2) =\cr
\sum_{r=0}^l\sum_{s=0}^m c_{l,m;r,s}^{(n-2)} &R_{l-r,m-s}^{(n-2+r+s)}
(X_2,X_2^{\ast},Q; q^2)
R_{r,s}^{(n-3)}(X_1,X_1^{\ast}, Q-X_2X_2^{\ast}; q^2)\cr
&\ot 
R_{l-r,m-s}^{(n-2+r+s)} (Y_2,Y_2^{\ast},D; q^2) Y_1^r (Y_1^{\ast})^s\cr}\tag 3.5.7
$$
in which $(\a=n-2)$
$$
\eqalign{
c_{l,m;r,s}^{(\a)} &= c_{l,m}^{(\a)} a_{l,m,r,s}^{-1} (c_{r,s}^{(\a-1)})^{-1}\cr
&= {{1-q^{2(\a+r+s+1)}}\over  {1-q^{2(\a+1)}}} {{c_{l,m}^{(\a)}}\over {c_{l-r,m-s}^{(\a+r+s)} 
c_{r,s}^{(\a-1)}}}\cr}\tag 3.5.8
$$
and $c_{l,m}^{(\a)}$ is as in $(3.3.2)$. Here we employed the following
\smallskip
{\bf Notation } : For $\a >0$ and $l,m\in\z+$ we put
$$
R_{l,m}^{(\a)}(A,B,C;q) =\left\{\aligned C^m A^{l-m} P_m^{(\a,l-m)}({{C-AB}
\over C};q)
\qquad (l\geq m)\\
C^l P_l^{(\a,m-l)}({{C-AB}\over C};q) B^{m-l}\qquad (l\leq m).\endaligned\right.\tag 3.5.9
$$
Then $R_{l,m}^{(\a)}(z_n,w_n,1;q^2) \equiv R_{l,m}^{(\a)}(z_n,w_n;q^2)$.
\medskip
Let us have a closer look at these polynomials $R_{l,m}^{(\a)}(A,B,C;q)$:
$$
R_{l,m}^{(\a)}(A,B,C;q) = \left\{\aligned C^m A^{l-m}\sum_{k=0}^m {{(q^{-m};q)_k
(q^{\a+l+1};q)_k}\over {(q^{\a +1};q)_k (q;q)_k}} (q {{C-AB}
\over C})^k\\
C^l \sum_{k=0}^l {{(q^{-l};q)_k(q^{\a+m+1};q)_k}\over {(q^{\a +1};q)_k 
(q;q)_k}} (q {{C-AB}\over C})^k B^{m-l}\endaligned\right.
$$
in the respective cases $l\geq m$ and $l\leq m$. So the polynomials are rational in 
$q^{\a}$. Hence (see also $(3.3.2)$) both sides of $(3.5.7)$ are rational functions of 
$q^{\a}$. Multiplying with a suitable factor, we will obtain from $(3.5.7)$ an identity 
which is polynomial in $q^{\a}$ and holds true for $\a = 1,2,\hdots$. But then obviously 
the identity is true for all $\a>0$.\hb 
Finally, let $\sigma : \Cal Y\to
\Cal Y$ be the automorphism that sends $Y_1$ to $-q Y_1^{\ast}$, 
$Y_1^{\ast}$ to $-q^{-1} Y_1$, and fixes $Y_2$ and 
$Y_2^{\ast}$
(which comes down to interchanging $\b$ and $\g$ in ${\Cal Z}_2$). 
If we now apply $id\ot \sigma$ to $(3.5.7)$ we end up with:
\smallskip
\proclaim\nofrills {\bf Theorem 3.5.8} : Suppose we are given the abstract complex 
$\ast$-algebras $\Cal X$ and $\Cal Y$ with generators $X_1,X_2,X_1^{\ast},
X_2^{\ast},Q$ respectively $Y_1,Y_2,Y_1^{\ast},
Y_2^{\ast},D$, relations $(3.5.6)$ and $\ast$-structures $(3.5.5)$.\hb
Then, for arbitrary $\a>0$ and arbitrary $l,m\in\z+$ we have the following 
addition formula for $q$-disk polynomials:
$$
\eqalign{
R_{l,m}^{(\a)}(-q X_1\ot &Y_1^{\ast}+X_2\ot Y_2, -qX_1^{\ast}\ot Y_1
+X_2^{\ast}\ot Y_2^{\ast}, QD;q^2) =\cr
\sum_{r=0}^l\sum_{s=0}^m c_{l,m;r,s}^{(\a)} &R_{l-r,m-s}^{(\a+r+s)}
(X_2,X_2^{\ast},Q; q^2)
R_{r,s}^{(\a-1)}(X_1,X_1^{\ast}, Q-X_2X_2^{\ast}; q^2)\cr
&\ot
(-q)^{r-s} R_{l-r,m-s}^{(\a+r+s)} (Y_2,Y_2^{\ast},D; q^2) Y_1^s (Y_1^{\ast})^r \cr}
\tag 3.5.10
$$
in which we used the notations $(3.5.8)$ and $(3.5.9)$. 
\endproclaim
\medskip
{\bf Remark } : For $\a=n-2$ this is in fact an identity in $\Z\ot {\Cal Z}_2$,
which we can rewrite as an identity in $\wi{\Z} \ot \wi{{\Cal Z}_2}$
by putting $Q=D=1$ in $(3.5.10)$.
For general $\a > 0$ we can do something similar (the relations among the 
generators are then given by $(3.5.6)$ but with $Q$ and $D$ equal to $1$).


\medskip

\Refs
\widestnumber\key{NYM}
\ref\key AA \by G.E. Andrews, R. Askey \paper Enumeration of partitions: The
role of Eulerian series and $q$-orthogonal polynomials \inbook in "Higher
combinatorics", M. Aigner (ed.), Reidel \yr 1977 \pages 3- 26\endref
\ref\key B \by G.M. Bergman \paper The diamond lemma for ring theory \jour Adv. in Math. 
\yr 1978 \vol 29 \pages 178-218\endref
\ref \key D \by M.S. Dijkhuizen \book On compact quantum groups and quantum homogeneous 
spaces  \publ PhD thesis University of Amsterdam \yr February 1994\endref
\ref\key DK \by M.S. Dijkhuizen, T.H. Koornwinder \paper CQG algebras: a direct algebraic
approach to compact quantum groups \jour Report AM-R9401, CWI, Amsterdam
\yr 1994 \moreref \jour hep-th/9406042; to appear in Lett. Math. Phys\endref
\ref\key Ko1 \manyby H.T. Koelink \paper On $\ast$-representations of the Hopf
$\ast$-algebra associated with the quantum group $U_q(n)$ \jour Compositio
Math. \yr 1991 \vol 77 \pages 199-231\endref
\ref \key Ko2 \bysame \paper The addition formula for continuous
$q$-Legendre polynomials and associated sphe- rical elements on the $SU(2)$
quantum group related to Akey-Wilson polynomials \jour SIAM J. Math. Anal.
\vol 25, No.1 \pages 197-217 \yr 1994\endref
\ref \key Ko3 \bysame \paper Askey-Wilson polynomials and the quantum 
$SU(2)$ group: survey and applications \jour to appear in Acta Appl. 
Math\endref
\ref \key Koo1 \manyby T.H. Koornwinder \paper The addition formula for 
Jacobi polynomials I, Summary of the results \jour Nederl. Akad. Wetensch.
Proc. Ser. A 75 \yr 1972 \pages 188-191\endref
\ref \key Koo2 \bysame \paper The addition formula for Jacobi polynomials, 
Parts II and III \jour Reports TW 133/72 and 135/72, Math. Centrum, Amsterdam 
\yr 1972\endref
\ref \key Koo3 \bysame \paper Orthogonal polynomials in
connection with quantum groups \inbook in "Orthogonal Polynomials: Theory
and Practice", P. Nevai (ed.) \publ NATO ASI Series, vol 294, Kluwer, 1990
\pages 257-292 \endref
\ref\key Koo4 \bysame \paper Positive convolution structures  
associated with quantum 
groups \inbook in "Probability Measures on Groups X", H. Heyer (ed.) \publ 
Plenum \yr 1991\pages 249-268\endref
\ref \key Koo5 \bysame \paper The addition formula for little
$q$-Legendre polynomials and the $SU(2)$ quantum group \jour SIAM J. Math.
Anal. \vol 21 \yr 1991 \pages 295-301\endref
\ref\key Koo6 \bysame \paper Askey-Wilson polynomials as zonal spherical functions
on the $SU(2)$ quantum group \jour SIAM J. Math. Anal. \vol 24 \yr 1993
\pages 795-813\endref
\ref\key Koo7 \bysame \paper General compact quantum groups, a 
tutorial \jour Report 94-06, Math. Preprint series, Dept. of Math. and 
Comp. Sci., University of Amsterdam (1994)\moreref \jour hep-th/ 9401114 
\moreref
\inbook part of the paper "Compact quantum groups and $q$-special
functions" in "Representations of Lie groups and quantum groups", 
V. Baldoni \& M. Picardello (eds.), Pitman Research Notes in Mathematics 
Series 311, Longman Scientific \& Technical, 1994.\endref
\ref\key N1 \manyby M. Noumi \paper Quantum groups and $q$-orthogonal polynomials.
Towards a realization of Askey -Wilson polynomials on $SU_q(2)$ \inbook
in "Special Functions", M. Kashiwara \& T. Miwa (eds.)\publ ICM-90
Satellite Conference Proceedings, Springer \yr 1991 \pages 260-288\endref
\ref \key N2 \bysame  \paper Macdonald's symmetric polynomials as zonal
spherical functions on some quantum homogeneous spaces  \jour to appear in 
Adv. in Math\endref
\ref\key NM \by M. Noumi, K. Mimachi \paper Askey-Wilson polynomials and
the quantum group $SU_q(2)$ \jour Proc. Japan Acad., Ser.A \vol 66 \yr 1990
\pages 146-149\endref
\ref\key NYM \by M. Noumi, H. Yamada, K. Mimachi \paper Finite dimensional
representations of the quantum group $GL_q(n, {\Bbb C})$ and the zonal 
spherical functions on $U_q(n-1)\backslash U_q(n)$ \jour Japan. J. Math.
(N.S.) \vol 19 \yr 1993 \issue 1 \pages 31-80\endref
\ref\key RTF \by N.Yu. Reshetikhin, L.A. Takhtadzhyan, L.D. Faddeev \paper Quantization of Lie 
groups and Lie algebras \jour Algebra and Analysis \vol 1 \yr1989 \pages 178-206 \moreref \jour 
Leningrad Math. Jnl. \vol 1 \yr 1990\pages 193-225\endref
\ref\key S \by R.L. $\breve{\text{S}}$apiro \paper Special functions connected with representations
of the group $SU(n)$ of class I relative to $SU(n-1)\ (n\geq 3)$ \jour Izv.
Vys$\breve{\text{s}}$. U$\breve{\text{c}}$ebn. Zaved. Mathematika \vol 71 
\yr 1968\pages 9-20 (Russian)
\moreref \paper AMS Translation Series 2 \vol 113 \yr 1979 \pages 
201-211 (English)\endref 
\ref\key Sw \by M.E. Sweedler \book Hopf algebras \publ W.A. Benjamin 
\publaddr New York \yr 1969\endref
\ref\key VS \by L.L. Vaksman, Ya. S. Soibel'man \paper The algebra of functions on the 
quantum group $SU(n+1)$ and odd-dimensional quantum spheres \jour Leningrad Math. Jnl. 
\vol 2 \yr 1991 \pages 1023-1042\endref
\bye